\documentclass[11pt]{article}
\usepackage{amsfonts}
\usepackage{amssymb}
\usepackage{amsbsy}
\usepackage{amsthm}

\newtheorem{thm}{Theorem}[section]
\newtheorem{cor}[thm]{Corollary}
\newtheorem{lem}[thm]{Lemma}
\newtheorem{prop}[thm]{Proposition}
\theoremstyle{definition}
\newtheorem{defn}[thm]{Definition}
\theoremstyle{remark}
\newtheorem{rem}[thm]{Remark}

\newcommand{\bt}{\begin{thm}}
\newcommand{\et}{\end{thm}}
\newcommand{\bc}{\begin{cor}}
\newcommand{\ec}{\end{cor}}
\newcommand{\bl}{\begin{lem}}
\newcommand{\el}{\end{lem}}
\newcommand{\bp}{\begin{prop}}
\newcommand{\ep}{\end{prop}}
\newcommand{\bd}{\begin{defn}}
\newcommand{\ed}{\end{defn}}
\newcommand{\br}{\begin{rem}}
\newcommand{\er}{\end{rem}}
\newcommand{\bpr}{\begin{proof}}
\newcommand{\epr}{\end{proof}}

\newcommand{\bi}{\begin{itemize}}
\newcommand{\ei}{\end{itemize}}
\newcommand{\be}{\begin{enumerate}}
\newcommand{\ee}{\end{enumerate}}
\newcommand{\ds}{\displaystyle}
\newcommand{\ba}{\begin{array}}
\newcommand{\ea}{\end{array}}
\newcommand{\beq}{\begin{equation}}
\newcommand{\eeq}{\end{equation}}
\newcommand{\beqa}{\begin{eqnarray}}
\newcommand{\eeqa}{\end{eqnarray}}

\newcommand{\N}{{\mathbb N}}

\newcommand{\R}{{\mathbb R}}
\newcommand{\C}{{\mathbb C}}
\newcommand{\T}{{\mathbb T}}
\newcommand{\PP}{{\mathbb P}}

\newcommand{\spec}{\mathrm{spec}}

\begin{document}

\title{\bf Matrix orthogonal polynomials whose derivatives are also orthogonal}
\author{M.J. Cantero $^a$, L. Moral $^b$, L. Vel\'azquez $^c$
\thanks{The work of the authors was supported, in part, by a research grant from the
Ministry of Education and Science of Spain, project code MTM2005-08648-C02-01,
and by Project E-64 of Diputaci\'on General de Arag\'on (Spain).} \\
\small{Departamento de Matem\'atica Aplicada, Universidad de Zaragoza, Spain} \\
\small{
(a) \texttt{mjcante@unizar.es} \
(b) \texttt{lmoral@unizar.es} \
(c) \texttt{velazque@unizar.es}}
}
\maketitle

\kern-10pt

\begin{abstract}

In this paper we prove some characterizations of the matrix orthogonal polynomials whose
derivatives are also orthogonal, which generalize other known ones in the scalar case. In
particular, we prove that the corresponding orthogonality matrix functional is
characterized by a Pearson-type equation with two matrix polynomials of degree not
greater than 2 and 1. The proofs are given for a general sequence of matrix orthogonal
polynomials, not necessarily associated with an hermitian functional. However, we give
several examples of non-diagonalizable positive definite weight matrices satisfying a
Pearson-type equation, which show that the previous results are non-trivial even in the
positive definite case.

A detailed analysis is made for the class of matrix functionals which satisfy a
Pearson-type equation whose polynomial of degree not greater than 2 is scalar. We
characterize the Pearson-type equations of this kind that yield a sequence of matrix
orthogonal polynomials, and we prove that these matrix orthogonal polynomials satisfy a
second order differential equation even in the non-hermitian case. Finally, we prove and
improve a conjecture of Dur\'{a}n and Gr\"{u}nbaum concerning the triviality of this class in the
positive definite case, while some examples show the non-triviality for hermitian
functionals which are not positive definite.

\end{abstract}


\medskip

\noindent{\it (2000) AMS Mathematics Subject Classification}: 42C05.

\vfill\eject

\section{Introduction}

\medskip

The results published by Dur\'{a}n in \cite{D97} can be considered the starting point for a
general study of matrix orthogonal polynomials satisfying differential equations. After
\cite{D97}, many other papers on the subject have appeared trying to find the
similarities and main differences with respect to the classical and semi-classical scalar
orthogonal polynomials (see \cite{CMV04,CMV05,CG05,CG,DG04,DG05c,DG05b,GI}). In spit of
these efforts, a complete Bochner-type classification of matrix orthogonal polynomials
satisfying second order differential equations similar to the scalar case (see
\cite{ABMP92,Boch29}) is far from being obtained.

However, the are many other differential properties that characterize the classical
scalar orthogonal polynomials and that could lead to interesting matrix generalizations.
These generalizations could clarify the structure of certain families of matrix
orthogonal polynomials, being a source of properties for such families, as in the scalar
case. Eventually, the understanding of these other differential properties could shed
light on the structure of some families of matrix orthogonal polynomials satisfying
differential equations, helping to find classification theorems.

It is well known that, apart from the second order differential equation, the classical
scalar orthogonal polynomials $(P_n)$ can be characterized by the orthogonality of their
derivatives $(P_{n+1}')$ (see \cite{BLN87,Ch78,HVR85,Mar91,Mar93}) or, equivalently, by a
linear relation between $P_n$ and $P_{n+1}',P_n',P_{n-1}'$ (see \cite{MBP94}). Also,
these properties are equivalent to a Pearson-type equation for the corresponding
orthogonality functional (see \cite{Ch78,Mar91,Mar93,Sho39}). The main objective of this
paper is to prove that the equivalence between these three properties hold in the matrix
case too (see Theorem \ref{MAIN}).

The proofs of the above equivalences are given for any sequence of matrix orthogonal
polynomials, not necessarily related to an hermitian weight matrix. Consequently, the
Pearson-type equation must involve a distributional derivative. The distributional
definition of the derivative not only permits to prove the results in a more general
context, but unifies many different situations that otherwise would require a separate
discussion. The reason is that the distributional Pearson-type equation takes care, not
only of the first order differential equation for the weight, but of the necessary
additional boundary conditions too (see Remark \ref{DISTR}). So, the introduction of the
distributional derivative becomes an advantage that permits to obtain more general
results and, at the same time, in a simpler and more elegant way.

Diagonalizable matrix orthogonal polynomials (we will be more precise about this concept
later) are nothing really different from scalar orthogonal polynomials. So, the relevance
of the results proved in this paper depends on the existence of non-diagonalizable
examples of matrix orthogonal polynomials whose derivatives are also orthogonal. Examples
2, 3 and 4 show that there are non-diagonalizable positive definite weight matrices whose
orthogonal polynomials enjoy such a property.

The weight matrix given in Example 2
$$
e^{-x^2} \pmatrix{1 + |a|^2x^2 & ax \cr \bar ax & 1} dx, \quad x\in\R, \quad
a\in\C\setminus\{0\},
$$
appeared previously in \cite{DG05b} as an archetype of positive definite weight matrices
whose orthogonal polynomials satisfy a second order differential equation. Curiously, the
authors declare in \cite{DG05b}, Section 7, Proposition 7.3, that the derivatives of
these matrix orthogonal polynomials are no longer orthogonal with respect to any weight
matrix, arguing that a contradiction appears when supposing a three term recurrence
relation for such derivatives. However, if one makes the proposed computations in
\cite{DG05b}, Proposition 7.3, no contradiction appears! Indeed, we will see that this
weight matrix satisfies a Pearson-type equation that, according to Theorem \ref{MAIN},
implies the orthogonality of the derivatives of its orthogonal polynomials. Even more, we
will find the positive definite weight matrix that gives the orthogonality of these
derivatives.

The purpose of \cite{DG05b}, Section 7, was to show that the equivalent characterizations
of the classical scalar orthogonal polynomials do not necessarily hold for matrix
orthogonal polynomials satisfying second order differential equations. It seems that the
authors were not too lucky in the choice of the weight matrix since, if they had chosen
the other example that they present, namely,
$$
e^{-x^2} \pmatrix{1 + |a|^2x^4 & ax^2 \cr \bar ax^2 & 1} dx, \quad x\in\R, \quad
a\in\C\setminus\{0\},
$$
they would have succeeded. The reason is that, as can be easily checked, this other
weight does not satisfy the required Pearson-type equation and, then, Theorem \ref{MAIN}
implies that the derivatives of its orthogonal polynomials can not be orthogonal.

A particular class of the family of matrix orthogonal polynomials with orthogonal
derivatives permits a deeper analysis. This is the class corresponding to a Pearson-type
equation involving a scalar polynomial $\alpha$ under the derivative. These matrix
orthogonal polynomials can be classified analogously to the classical scalar case,
according to the roots of $\alpha$: Hermite (no roots), Laguerre (a simple root), Jacobi
(two different roots) or Bessel-type (a double root). Moreover, a change of variable can
reduce the different types to the canonical cases $\alpha(x)=1,x,1-x^2,x^2.$

For this special class we develop explicit formulas for the related matrix parameters,
such as the norm of the monic orthogonal polynomials, the coefficients of the three term
recurrence relation or the coefficients of the linear relation between the polynomials
and their derivatives. These formulas, although generalizations of the known ones in the
classical scalar case, are more intricate due to the non-commutativity of the matrix
product. However, they are very useful since they allow to characterize the Pearson-type
equations that have a quasi-definite solution. In other words, if a matrix functional
satisfies this kind of Pearson-type equation, we have a criterion to know if it generates
a sequence of orthogonal polynomials (see Theorem \ref{PTEZCL}). Notice that the
importance of this result relies on the fact that we are dealing with general matrix
functionals and not only with positive definite weight matrices, since the last ones
always have an associated sequence of matrix orthogonal polynomials.

We also prove that the matrix orthogonal polynomials of the above class satisfy a second
order differential equation with polynomial coefficients (see Theorems \ref{DifEq} and
\ref{DifEq2}). The result is again true no matter if the corresponding orthogonality
matrix functional is hermitian or not. This is one of the novelties of this result, since
the previous works on differential equations for matrix orthogonal polynomials always
dealt with the hermitian case only. Indeed, if we believe a conjecture formulated by
Dur\'{a}n and Gr\"{u}nbaum in \cite{DG05}, this discovering is only relevant for the functionals
of the referred class that are not positive definite. This conjecture says that every
positive definite weight matrix in this class is diagonalizable. We present a proof of
this conjecture (see Corollary \ref{DGconj}).

The above conjecture was supported on a partial proof given in \cite{DG05}, that was
incomplete due to the strong assumptions made there. First of all, it was supposed that
the coefficients of the matrix polynomial appearing in the Pearson-type equation commute.
Second, the proof was given separately for each of the canonical types of hermitian
weight matrices that in the scalar case are positive definite: $\alpha(x)=1,x,1-x^2.$ So,
the case $\alpha(x)=x^2$ is not considered, although the authors do not prove its
incompatibility with a positive definite weight in the matrix case too. Finally, there is
another less evident inconvenient. If $\alpha$ has a complex root, the required change of
variable to arrive at a canonical situation destroys in general the hermiticity of the
weight matrix. This means that, apart from the previous restrictions, the proof is only
valid for the case of $\alpha$ with real roots. Our proof avoid all these problems. Even
more, we get a result that improves the one conjectured in \cite{DG05} (see Theorem
\ref{DGC3}). In spite of this result, the non-triviality of the class under consideration
is ensured by the existence of non-diagonalizable matrix orthogonal polynomials in such a
class, even in the hermitian case (see \cite{CMV05,DG05} and Example 5 of this paper).

The exposition of the above results will be structured in the following way along the
paper. Section 2 introduces the notation, as well as some preliminary results and
considerations that will of interest for the rest of the paper. In Section 3 we study the
matrix orthogonal polynomials $(P_n)$ with respect to a functional satisfying a
Pearson-type equation with two matrix polynomials of degree not greater than 2 and 1. We
prove that such a Pearson-type equation is equivalent to the orthogonality of the
derivatives $(P_{n+1}')$ and, also, to a linear relation between $P_n$ and
$P_{n+1}',P_n',P_{n-1}'.$ Some two-dimensional non-diagonalizable examples of positive
definite weight matrices whose orthogonal polynomials satisfy these properties are
presented at the end of the section. Section 4 is devoted to the analysis of the special
case in which the polynomial under the derivative in the Pearson-type equation is a
scalar one. We obtain the characterization of the Pearson-type equations of this kind
with quasi-definite solutions, the differential equation for the related matrix
orthogonal polynomials and the proof of the Dur\'{a}n-Gr\"{u}nbaum conjecture, finishing with
some non-diagonalizable examples. Finally, in Section 5 we discuss the relation of the
above results with other ones in the literature about second order differential equations
for matrix orthogonal polynomials.

\smallskip

\section{The Basics}

\medskip

We start with some notations and a summary of basic results that we will use in the rest
of the paper.

\medskip

In what follows, $\C^m$ will be the set of complex vectors of $m$ components and
$\C^{(m,m)}$ the set of $m \times m$ complex matrices. We shall denote by $\PP^{(m)}$ the
$ \C^{(m,m)}$-left-module
$$
\PP^{(m)}=\left\{\sum_{k=0}^n \alpha_k x^k \,\bigg|\, \alpha_k\in \C^{(m,m)}, \; n\in\N
\right\},
$$
and by means of $\PP^{{(m)}'}$ the $\C^{(m,m)}$-right-mo\-du\-le Hom$\left(\PP^{(m)},
\C^{(m,m)}\right).$ $\PP_n^{(m)} $ will be the subset of matrix polynomials  of
$\PP^{(m)}$ with degree not greater than $n.$ In the scalar case ($m=1$) we will just write
$\PP^{(1)}=\PP$ and $\PP^{(1)}_n=\PP_n$.

For all $P\in \PP^{(m)}$ and  $u\in \PP^{{(m)}'}$ the duality bracket is defined by
$\langle P, u  \rangle = u\left(P\right)$ and it verifies the usual bilinear properties.

For $k\in \N$ and $u\in\PP^{{(m)}'}$ the linear functional $u x^k I\in \PP^{{(m)}'}$ is
given by
$$
\langle P, u x^kI \rangle = \langle x^k P, u \rangle,
$$
where $I$ denotes the $m\times m$ identity matrix. A linear extension gives the
right-product $u Q \in \PP^{{(m)}'}$ for $u\in \PP^{{(m)}'},$ $Q\in \PP^{(m)},$ with
$Q(x)=\sum_{k=0}^n q_k x^k,$ $q_k\in\C^{(m,m)}$, in the following way:
$$
\langle P, u Q \rangle =  \sum_{k=0}^n  \langle x^k P, u  \rangle q_k.
$$
Similarly, the left-product $Q u \in \PP^{{(m)}'}$ is defined by
$$
\langle P,Qu \bigr> =  \langle PQ,u \rangle.
$$

Every functional $u \in \PP^{{(m)}'}$ induces a matrix inner product in
$\mathbb{P}^{(m)}$ given by $ \langle P, Q  \rangle_u =  \langle P, u Q^* \rangle, $
where $Q^*(x)=\sum_{k=0}^n q_k^* x^k$ and $q_k^*$ is the adjoint matrix of $q_k$.
This matrix inner product enjoys the standard sesquilinear properties. The orthogonality
with respect to $u$ means the orthogonality with respect to this inner product.

The functional $u^*$ is defined by
$$
 \langle P, u^* Q  \rangle  = \langle Q^*, u P^* \rangle^*,
$$
and we will say that $u$ is an hermitian functional if $u=u^*.$ In  this case $\langle P,
uP^*\rangle$ is hermitian for any $P\in \PP^{(m)}.$ We will say that an hermitian
functional $u$ is positive definite if $\langle P, uP^* \rangle$ is positive definite for
every $P\in\mathbb{P}^{(m)}$ with $\det P \neq 0. $ In what follows we denote this
condition by $u>0.$ In the same way, for a positive definite matrix $A$ we will write
$A>0.$

\medskip

We denote by $\mu_k = \langle x^kI, u \rangle$ the $k$-th moment with respect to
$u\in\PP^{{(m)}'}$. Given a sequence $\left(\mu_k\right)_{k \geq 0}$ in $\C^{(m,m)},$
there exists a unique  $u\in\mathbb{P}^{(m)'}$ such that $ \langle x^kI, u \rangle  =
\mu_k.$

If $u\in \PP^{{(m)}'}$ has moments $\left(\mu_k\right)_{k\geq0}, $ we say that $u$ is
quasi-definite (or non-singular) if $\det \Delta_n \not=0$ for $n\geq 0,$ where
$\Delta_n$ is the Hankel-block matrix
$$
\Delta_n = \pmatrix{\mu_0 &\mu_1&\dots &\mu_n \cr
                    \mu_1 &\mu_2&\dots &\mu_{n+1} \cr
                    \dots&\dots&\dots&\dots  \cr
                    \mu_n&\mu_{n+1}&\dots&\mu_{2n}\cr}.
$$
Notice that $u$ is hermitian if and only if $\mu_n = \mu_n^*$ for $n\geq0$, or,
equivalently, $\Delta_n= \Delta_n^*$ for $n\geq0.$

The interest of the quasi-definite functionals relies on the following result (see
\cite{D95,DVAss95,San92}).

\begin{thm}\label{S1T1}
$u\in\mathbb{P}^{(m)'}$  is quasi-definite if and only if there exists a sequence
$(P_n)_{n\geq 0}$ of left  orthogonal  matrix polynomials with respect to $u,$ that is:

\item{(i)} $P_n \in\mathbb{P}^{(m)},\ \deg P_n = n.$
\item{(ii)} The leading coefficient of $P_n$ is non-singular.
\item{(iii)} $ \langle x^k P_n, u \rangle = E_n\delta_{nk},$ with $E_n$ non-singular, for $0 \leq k \leq n$.
\smallskip

\noindent Moreover, the sequence $(P_n)_{n\geq 0}$ is unique up to non-singular left
matrix factors and verifies a recurrence  relation
$$
x P_n(x) = \alpha_n P_{n+1}(x) + \beta_n P_n(x) + \gamma_n  P_{n-1}(x),
$$
where $P_0\in \mathbb{C}^{(m,m)}$ is non-singular, $P_{-1}= 0$ and $\alpha_n,$ $\beta_n,$
$\gamma_n \in \mathbb{C}^{(m,m)},$ with $\alpha_n,$ $\gamma_n$ non-singular.
\end{thm}

The last result of this theorem has a converse (Favard's Theorem): for any sequence
$(P_n)_{n\geq 0}$ verifying the above recurrence relation there exists a unique (up to
non-singular right matrix factors) quasi-definite functional $u$ such that
$(P_n)_{n\geq0}$ is its sequence of left orthogonal  matrix polynomials (see
\cite{D95,DVAss95,San92}). Analogously we can define the right orthogonal matrix
polynomials with respect to $u$, which are the adjoints of the left orthogonal
polynomials associated with $u^*$. In what follows we will consider only left orthogonal
matrix polynomials, and we will call them just matrix orthogonal polynomials (MOP).

\br\label{ATH}
Given a functional $u\in\mathbb{P}^{(m)'}$, we can normalize the corresponding MOP by
choosing the only monic ones $(P_n)_{n\geq0}$. In what follows we will assume this
choice, so, a unique sequence of non-singular matrices $(E_n)_{n \geq 0}$, $E_n = \langle
x^nP_n, u \rangle$, is associated with any quasi-definite functional $u$. Also, $\beta_n$
and $\gamma_n$ will denote the matrix coefficients of the related recurrence relation
$$
xP_n(x) = P_{n+1}(x) + \beta_n P_n(x) + \gamma_n P_{n-1}(x).
$$

Similarly, given a sequence MOP, we can normalize the corresponding functional $u$ in
different ways, for instance, by requiring $\langle I, u \rangle = I$. However, we will
not fix the normalization for the moment because the  most convenient one depends on the
problem that we wont to study.

\er

\medskip

In the case of non quasi-definite functionals, the full sequence of MOP does not exist.
Nevertheless, we have the following general result.

\bp \label{S1P1}
For every $u\in \PP^{{(m)}'}$ the following statements are equivalent:

\item{(i)} $\Delta_0,\dots,\Delta_n$ are non-singular.

\item{(ii)} There exists a  finite segment $\bigl(P_k\bigr)_{k=0}^n$ of monic MOP with
respect to $u,$ ${}\kern15pt$ that is:
\smallskip

(a) $P_k \in\mathbb{P}^{(m)}, \ \deg P_k = k.$
\smallskip

(b) $\langle x^j P_k, u \rangle = E_k\delta_{kj},$ with $E_k$ is non-singular, for
$0 \leq j \leq k \leq n.$
\smallskip

\noindent Moreover, under the above conditions, the segment $\bigl(P_k\bigr)_{k=0}^n$ is
unique and  there exists a unique monic polynomial $P_{n+1}$ whit $\deg P_{n+1}= n+1$
such that $\langle x^j P_{n+1}, u \rangle = 0$ for $0 \leq j \leq n$.
\ep

\bpr
Suppose that $\Delta_0,\dots,\Delta_n$ are non-singular. If $P_k(x) =\sum_{i=0}^k
\pi_i^{(k)} x^i,$ $\pi_i^{(k)} \in \mathbb{C}^{(m,m)},$ then, $\langle x^j P_k, u \rangle
= \sum_{i=0}^k \pi_i^{(k)} \mu_{i+j}.$ Choosing $\pi_k^{(k)}=I,$ the system $\sum_{i=0}^k
\pi_i^{(k)} \mu_{i+j}= 0, \;j=0,\dots, k-1, $ can be represented as
$$
\pmatrix{\pi_0^{(k)},&\kern-7pt\pi_1^{(k)},&\kern-7pt\dots, &\kern-7pt\pi_{k-1}^{(k)}\cr}
\Delta_{k-1} = - \pmatrix{\mu_k,& \kern-7pt \mu_{k+1},&\kern-7pt\dots,& \kern-7pt\mu_{2k-1}\cr},
$$
which has a unique solution for $k=0,1,\dots,n+1.$

On the other  hand, $E_k$ is non singular for $k=0,1,\dots,n.$ In fact, we have
$
\langle x^j P_k, u \rangle = E_k \delta_{kj}, \;j=0,\dots,k,\; k=0, \dots, n,
$
and, so,
$$
\pmatrix{
\pi_0^{(k)},&\kern-7pt\pi_1^{(k)},&\kern-7pt\dots,&\kern-7pt\pi_{k-1}^{(k)},&\kern-7ptI
}
\Delta_k = \pmatrix{ 0,&\kern-7pt 0,&\kern-7pt\dots,&\kern-7pt0,&\kern-7ptE_k }.
$$
If $E_k$ is singular, there exists $v\in\mathbb{C}^{m}\setminus\{0\}$ such that
$v^TE_k=0.$ Hence,
$$
\pmatrix{
v^T\pi_0^{(k)},&\kern-7pt v^T\pi_1^{(k)},&\kern-7pt\dots,
&\kern-7pt v^T\pi_{k-1}^{(k)},&\kern-7pt v^T
}
\Delta_k =
\pmatrix{ 0,&\kern-7pt 0,&\kern-7pt\dots,&\kern-7pt0,&\kern-7pt0 },
$$
and this result contradicts the non-singularity of $\Delta_k$ for $k=0,\dots,n.$

For the converse, let us suppose that there exists a finite segment
$\bigl(P_k\bigr)_{k=0}^n$ of MOP with respect to $u$ with $E_k = \langle x^k P_k, u
\rangle.$ It is easy to see that the conditions $\langle x^j Q_k, u \rangle = E_k
\delta_{kj},$ $j=0,\dots k,$ where $Q_k \in\mathbb{P}^{(m)}_k,$ ensures that $Q_k=P_k,$
$\;k=0,\dots,n.$ Writing $ Q_k(x)=\sum_{i=0}^k \pi_i^{(k)} x^i,$ the above assertion
means that, for $k=0, \dots, n,$ the system
$$
\pmatrix{
\pi_0^{(k)},&\kern-7pt\pi_1^{(k)},&\kern-7pt\dots,&\kern-7pt\pi_{k-1}^{(k)},&\kern-7pt\pi_k^{(k)}
}
\Delta_k =
\pmatrix{ 0,&\kern-7pt 0,&\kern-7pt\dots,&\kern-7pt0,&\kern-7ptE_k }
$$
has a unique solution and, hence, $\Delta_k$ is non-singular.
\epr

Concerning the partial hermiticity of a functional, we have the following immediate result.

\bp\label{S1P2}
Let $u\in \PP^{{(m)}'}.$ If $(p_k)_{k=0}^n$ is a basis of $\PP_n^{(m)},$
$\Delta_n=\Delta_n^*$ if and only if $( \langle p_k, up_j^* \rangle )_{k,j=0}^n$ is
hermitian.

In particular, if $u$ has  a finite segment $(P_k)_{k=0}^n$ of MOP,
$$
\Delta_n = \Delta_n^* \;\; \Longleftrightarrow \;\; \langle P_k, uP_j^*\rangle =
E_k \delta_{kj}, \;\; E_k = E_k^*, \;\; 0 \leq j,k \leq n.
$$
\ep

\smallskip

The second assertion of the above proposition says that, when $\Delta_0,\dots,\Delta_n$
are non-singular, the condition $\Delta_n=\Delta_n^*$  means that the finite segments of
left and right orthogonal matrix polynomials are each one the hermitian adjoint of the
other one.

Also, for the hermitian positive definite functionals on $\PP^{(m)}_n$ we have the
following characterization.

\bp\label{S1P3}
Let $u\in \PP^{{(m)}'}.$ If  $(p_k)_{k=0}^n$ is a basis of $\PP^{(m)}_n,$ the following
statements are equivalent:

\item{(i)} $\Delta_n > 0.$

\item{(ii)} $(\langle p_k,up_j^* \rangle)_{k,j=0}^n >0.$

\item{(iii)} $u$ has a finite segment $(P_k)_{k=0}^n$ of MOP such that
$\langle P_k, uP_j^*\rangle = E_k \delta_{kj}$ ${}\kern18pt$ with $E_k > 0$ for
$0 \leq j,k \leq n.$

\item{(iv)} $ \langle P, uP^*\rangle > 0$ for any $P\in \PP_n^{(m)}$ such that $\det P\neq 0.$
\ep

\bpr
We only prove {\it (i) $\Leftrightarrow$ (iv)}, since the rest of equivalences are
immediate. For any matrix polynomial $P(x)=\sum_{i=0}^k A_i x^i,$
$A_i\in\mathbb{C}^{(m,m)},$ $k\leq n,$
$$
\langle P, uP^*\rangle =
\pmatrix{A_1 \, \dots \, A_k\cr} \Delta_k \pmatrix{A_1^*\cr\vdots\cr A_k^*\cr}.
$$
So, $\langle P, uP^*\rangle$ is hermitian if $\Delta_n$ is hermitian.
If $v\in \mathbb{C}^{m},$
\beq\label{VV}
v^* \langle P, uP^* \rangle v =
\pmatrix{v_0^* \, \dots \, v_k^*\cr} \Delta_k \pmatrix{v_0\cr\vdots\cr v_k\cr},
\qquad  v_i=A_i^* v.
\eeq
Then, if  $v\neq 0,$ $\det P\neq 0$ implies $v_i \neq 0$ for some $i.$ So, equality
(\ref{VV}) gives $v^* \langle P, uP^*\rangle v > 0$  if $\Delta_n>0.$

For the converse, if $\langle P, uP^*\rangle$ is hermitian for $P\in\mathbb{P}_n^{(m)}$
with $\det P\neq 0,$ $  \mu_{2k} = \langle x^kI, u x^k I \rangle = \mu_{2k}^*$ for $k\leq
n.$ Besides, $\mu_{2k-1}=\mu_{2k-1}^*$ for $k\leq n$ too, due to the identity
$\langle(x^k+x^{k-1})I,u(x^k+x^{k-1})I\rangle = \mu_{2k}+\mu_{2k-2}+2\mu_{2k-1}$.
Therefore $\Delta_n=\Delta_n^*.$

Suppose $\langle P, uP^* \rangle > 0$ for any $P\in\mathbb{P}_n^{(m)}$ with $\det P \neq
0.$ Let $\pmatrix{v_0 \dots v_k}, $ $v_i \in \mathbb{C}^{m},$ with $v_k\neq 0$ and $k\leq
n.$ We can always find $A_i \in \mathbb{C}^{(m,m)}$ such that $A_i^*v_k = v_i,$ $A_k=I.$
The polynomial $P(x) =\sum_{i=0}^k A_ix^i$ lies on $\PP^{(m)}_n$ and $\det P\neq 0.$ So,
relation (\ref{VV}) gives
$$
\pmatrix{v_0^* \, \dots \, v_k^*} \Delta_k \pmatrix{v_0\cr\vdots\cr v_k\cr} >0,
\quad \hbox{if } v_k\neq0, \quad k \leq n.
$$
This proves by induction that $\Delta_n >0.$
\epr

\br\label{POS}
Notice that, if $u$ is an hermitian and positive definite functional, then it is
quasi-definite. So, there exits the corresponding  sequence $(P_n)_{n\geq 0}$ of MOP with
$E_n$ hermitian and positive definite.

Similarly to the scalar case, the positive definite matrix functionals are those ones given by
\beq\label{SD}
\langle P, u \rangle =  \int P(x) \, dM(x),
\eeq
where $dM$ is a positive definite weight matrix on $\R$, that is, a positive definite
matrix of measures supported on the real line ($M(S)$ is positive semidefinite for any
Borel set $S\subset\R$) with finite moments $\int x^ndM(x),$ $n\geq0,$ and such that
$\int P(x)\,dM(x)\,P(x)^*$ is non-singular if $\det P\neq0$ (see \cite{D95}). This is,
for instance, the case of an absolutely continuous matrix of measures $dM(x)=W(x)\,dx$
with finite moments, $W(x)$ being semidefinite positive for any $x\in\R$ and non-singular
for infinitely many points of the real line.

In what follows we will identify any $m \times m$ matrix $dM$ of measures on $\C$ with
finite moments (not necessarily hermitian), and the functional $u\in\PP^{{(m)}'}$ defined
by (\ref{SD}). Thus, we will write $u=dM$ for such a functional.
\er

\medskip

A specially interesting family of matrix functionals is given by the functionals  which
satisfy a differential equation of Pearson-type (see \cite{CMV04,CMV05}). The definition
of this family requires the introduction of the derivative operator in the space
$\PP^{{(m)}'},$ which is the linear operator $D \colon \PP^{{(m)}'} \to \PP^{{(m)}'}$
such that
$$
\langle P, Du \rangle = - \langle P', u \rangle.
$$
The equality $D(u \Phi)  = (D u) \Phi + u \Phi'$ holds for all $u\in\PP^{(m)'}$ and
$\Phi\in\PP^{(m)}.$

\begin{defn}
Let $u\in \PP^{{(m)}'}.$ We say that $u \in \mathcal{P}$ or, equivalently, $u$ is a
$\mathcal{P}$-functional, if there exist $\Phi, \Psi\in \PP^{(m)},$ with $\det\Phi\not=0$,
such that
$$
D\left(u\Phi\right) = u\Psi \qquad \hbox{(Pearson-type equation)}
$$
If $\deg \Phi\leq p$ and $\deg \Psi \leq q,$ we say that $u\in \mathcal{P}_{p,q}$ or $u$
is a $\mathcal{P}_{p,q}$-functional. In both  cases we also say  that the corresponding
sequence of MOP belongs to the family  $\mathcal{P}$ or $\mathcal{P}_{p,q}$ respectively.
\end{defn}

\br \label{DET}
The condition $\det\Phi\neq0$ is imposed to avoid any triviality of the definition, ensuring that
it involves all the components $u_{ij} \colon \PP^{(m)} \to \C$ of $u=(u_{ij})_{i,j=0}^m.$
Notice that
$$
\det\Phi=0 \kern7pt \Longleftrightarrow \kern7pt \Phi v = 0
\,\hbox{ for some }\, v\in\C^m[x]\setminus\{0\}.
$$
In fact, if $\Phi v = 0$ for some $v\in\C^m[x]\setminus\{0\},$ then $0 = (\hbox{\rm
adj}\,\Phi) \Phi v = (\det \Phi) v$. To see the converse, remember that every
$\Phi\in\PP^{(m)}$ can be factorized as $\Phi = P \hat\Phi Q$, with
$\hat\Phi\in\PP^{(m)}$ diagonal and $P,Q\in\PP^{(m)}$ invertible, that is, $\det P, \det
Q \in \C\setminus\{0\}.$ Therefore, $\det\Phi=0$ implies $\det\hat\Phi=0$ and, since
$\hat\Phi$ is diagonal, $\hat\Phi v_0=0$ for some $v_0\in\C^m\setminus\{0\},$ which gives
$\Phi v = 0$ with $v=Q^{-1}v_0\in\C^m[x]\setminus\{0\}.$
\er

\vskip0pt

\br \label{DISTR} The distributional definition of the derivative operator $D$ implies
that, in general, the Pearson-type equation involves, not only a relation between
standard derivatives, but a boundary condition too. Consider, for instance, a functional
$u=W(x)\,dx,$ $x\in\Gamma,$ with $W$ an analytic matrix function on a regular curve
$\Gamma$ of the complex plane. Then, $Du=W'(x)\,dx+W(x)(\delta(x-a)-\delta(x-b))\,dx,$
where $a$ and $b$ are the initial and end points of $\Gamma$ respectively. So, if the
curve is open, together with the equality $(W\Phi)'=W\Psi,$ we need the boundary
condition $(W\Phi)(a)=(W\Phi)(b)= 0$ to ensure the Pearson-type equation $D(u\Phi) =
u\Psi.$ The case of a closed curve does not need an additional boundary condition since
we suppose that $W$ is analytic on $\Gamma.$ Moreover, in this case, the Pearson-type
equation holds even if $(W\Phi)' \neq W\Psi$ but $(W\Phi)'-W\Psi$ is analytic on the
region enclosed by $\Gamma,$ due to Cauchy's Theorem. The Pearson-type equation can be
satisfied if $W$ is only analytic on $\Gamma\setminus\{a,b\}$ but the limits
$(W\Phi)(a^+):=\lim_{t \to t_0}(W\Phi)(\gamma(t))$, $(W\Phi)(b^-):=\lim_{t \to
t_1}(W\Phi)(\gamma(t))$ exist, where $\gamma\colon[t_0,t_1]\to\Gamma$ is a
parametrization of $\Gamma,$ $a=\gamma(t_0),$ $b=\gamma(t_1)$. Then,
$$
D(u\Phi)=(W\Phi)'(x)\,dx+(W\Phi)(a^+)\,\delta(x-a)\,dx-(W\Phi)(b^-)\,\delta(x-b)\,dx,
$$
so, we get the Pearson-type equation adding to $(W\Phi)'=W\Psi$ the boundary conditions
$$
\ba{l}
(W\Phi)(a^+) = (W\Phi)(b^-)
\kern42pt \hbox{\rm closed curve},
\smallskip
\cr
(W\Phi)(a^+) = (W\Phi)(b^-) = 0
\qquad \hbox{\rm open curve}.
\ea
$$
The distributional derivative not only unifies all these cases, but allows to consider
more general situations, such as functionals defined by matrix measures supported on an
arbitrary subset of the complex plane. \er

\medskip

If $u\in \PP^{{(m)}'}$ is a $\mathcal{P}$-functional with a Pearson-type equation
$D\left(u\Phi\right) = u\Psi$, then, for every $\Omega\in\PP^{(m)}$,
\beq \label{PDE}
D \left( u \Phi \Omega \right) = u \left(  \Phi\Omega' + \Psi\Omega\right).
\eeq
Therefore, the set
$$
{\cal M}(u) = \{ \Phi\in \PP^{(m)} \mid D(u\Phi) = u\Psi, \ \Psi\in \PP^{(m)} \}
$$
is a right-ideal of $\PP^{(m)}$, but it is not necessarily principal, because the
euclidean division algorithm is not valid in $\PP^{(m)}.$ This is an obstacle to find
a canonical representative of ${\cal M}(u)$ that could lead to a classification of
$\mathcal{P}$-functionals similarly to the scalar case.

Notice that $\mathcal{P} = \bigcup_{p,q\geq 0} \mathcal{P}_{p,q},$ and
$\mathcal{P}_{p,q}\subset\mathcal{P}_{p',q'}$ if $p \leq p'$ and $q \leq q'$. The set
$$
{\cal M}_{p,q}(u) =
\{ \Phi\in \mathbb{P}_p^{(m)} \mid D(u\Phi)=u\Psi, \ \Psi\in\PP_q^{(m)}\}
$$
is not an ideal of $  \PP^{(m)},$ but a $\mathbb{C}^{(m,m)}$-right-submodule of
$\PP_p^{(m)}.$ Although it is finitely generated, it is not cyclic in general, what means
again a problem for finding a canonical representative of ${\cal M}_{p,q}(u).$

\medskip

\noindent {\bf Example 1.} Let us consider $u\in \PP^{(2)'}$ given by
$$
u = (1-x^2)\pmatrix{1+3x^2 & 2x \cr 2x & 1}dx, \quad x\in(-1,1).
$$
A direct computation shows that $u$ is a $\mathcal{P}_{3,2}$-functional with
$$
{\cal M}_{3,2}(u)=
{\rm span}_{\mathbb{C}^{(2,2)}}
\biggl\{ (1-x^2)I , x(1-x^2)\pmatrix{0&0\cr0&1} \biggr\}
$$
generated by two elements. Indeed, if
$$
\Phi(x) = (1-x^2) \Lambda_1 + x(1-x^2)\pmatrix{0&0\cr0&1} \Lambda_2,
\quad \Lambda_i\in \mathbb{C}^{(2,2)},
$$
then $D(u\Phi) = u\Psi$ with
$$
\Psi(x) =
\pmatrix{-2x & 2 \cr 2-6x^2 & -8x} \Lambda_1 +
\pmatrix{0 & 2x \cr 0 & 1-9x^2} \Lambda_2.
$$

We can get cyclic modules for $u$ by going down in the net
$(\mathcal{P}_{p,q})_{p,q\geq0},$ but there are two different ways to do it. From
the previous result we obtain
\bi

\item $u\in \mathcal{P}_{2,2}$ with
${\cal M}_{2,2}(u) = {\rm span}_{\mathbb{C}^{(2,2)}}\biggl\{(1-x^2)I\biggr\}$.

\item $u\in \mathcal{P}_{3,1}$ with
${\cal M}_{3,1}(u) =
{\rm span}_{\mathbb{C}^{(2,2)}}\biggl\{(1-x^2)\pmatrix{3&0\cr-2x&1} \biggr\}$.

\ei
In fact,
$$
\ba{l}
D\left(u(1-x^2)I\right) = u\pmatrix{-2x & 2 \cr 2-6x^2 & -8x},
\medskip
\cr
D\left(u(1-x^2)\pmatrix{3&0\cr-2x&1\cr}\right) = u \pmatrix{-10x&2\cr4&-8x\cr}.
\ea
$$
This splitting shows clearly the problem of classification of $\mathcal{P}$-functionals.
Moreover, we can not go down more than this in the net $(\mathcal{P}_{p,q})_{p,q \geq 0}$
since
$$
\ba{l}
{\cal M}_{2,1}(u)= {\cal M}_{2,2}(u) \cap {\cal M}_{3,1}(u) =
{\rm span}_{\mathbb{C}^{(2,2)}}\left\{(1-x^2)\pmatrix{0&0\cr0&1\cr}\right\},
\cr
{\cal M}_{1,2}(u)= {\cal M}_{3,0}(u) = {\cal M}_{0,3}(u) = \{0\},
\ea
$$
and, hence, $u \not\in \mathcal{P}_{p,q}$ for $p+q\leq 3.$

\medskip

Notice that the above problems of classification happen even for quasi-definite
functionals since our example was positive definite. However, if we restrict our
attention to quasi-definite functionals, there is a singular situation. As we will prove
later (see Theorem \ref{MC}), if $\Delta_0, \Delta_1, \Delta_2$ are non-singular for some
$u\in \mathcal{P}_{2,1},$ then ${\cal M}_{2,1}(u)$ is cyclic. This implies that we can
associate with each sequence of MOP in the family $\mathcal{P}_{2,1}$ a canonical
representative: the unique (up to non-singular right matrix factors) generator of ${\cal
M}_{2,1}(u),$ $u$ being the related orthogonality matrix functional.

\medskip

A way to solve the problem of classification of $\mathcal{P}$-functionals uses the fact
that ${\cal M}(u)$ always has a non-trivial scalar representative. In fact, choosing
$\Omega=\hbox{\rm adj}\,\Phi$ in (\ref{PDE}) gives $\Phi \Omega  = (\det \Phi)I,$ which
yields the following characterization (see \cite{CMV04,CMV05}).

\bp
The functional  $u\in\PP^{{(m)}'}$ belongs to the family $\mathcal{P}$ if and only if
there exist $\alpha \in \PP \setminus \{0\}$ and $\Psi\in \PP^{(m)}$ such that
$$
D (u \alpha I) = u \Psi.
$$
\ep

Notice that the set
$$
\widetilde{{\cal M}}(u) =
\{ \alpha \in \PP  \mid    D (u\alpha I) = u \Psi, \ \Psi\in\PP^{(m)} \}
$$
is a non-trivial bilateral ideal of $\PP$, which is, therefore, principal. So,
there exists  an $\alpha \in \PP\setminus\{0\}$, unique up non-trivial factors in $\C,$
that is generator of $\widetilde{{\cal M}}(u).$ This scalar generator can be used to
classify the $\mathcal{P}$-functionals.

\bd
Let $u\in\PP^{(m)'}$ be a $\mathcal{P}$-functional and let $\alpha\in\PP\setminus\{0\}$
be a generator of $\widetilde{{\cal M}}(u). $ The  class of $u$ is
$s=\max\{\deg\alpha-2,\deg\Psi-1\},$ where $\Psi\in\PP^{(m)}$ is such that
$D(u\alpha I)=u\Psi$.
\ed

The interesting $\mathcal{P}$-functionals are those ones that have a sequence of MOP,
that is, the quasi-definite $\mathcal{P}$-functionals. These  are called semi-classical
functionals (see \cite{CMV04,CMV05}). As in the scalar case, the semi-classical
functionals can be characterized by several differential properties of the corresponding
MOP.

\bt\label{CHSCL}
Let  $u\in \PP^{{(m)}'}$ be  quasi-definite and let $\left(P_n\right)_{n\geq 0}$  be
the associated sequence of MOP. Then, the following statements are equivalent:

\item{(i)}    $u \in \mathcal{P}.$

\item{(ii)}   There exist $\alpha \in \PP\setminus\{0\}$ and
              $\Theta_j^{(n)} \in \C^{(m,m)}$ such that
$$
\kern25pt
\alpha(x) P'_{n+1}(x)  = \sum_{j=-s}^{\deg\alpha} \Theta_j^{(n)} P_{n+j}(x)
\qquad  \hbox{(structure  relation)}
$$
\kern22pt with $s\geq \max\{\deg\alpha-2,0\}$ independent of $n$ and
$\Theta_{-s}^{(n)} \not= 0$ for some
${}\kern22pt n \geq s.$

\item{(iii)}  There exist
              $a\in\PP\setminus\{0\}$, $b\in\PP$ and
              $\Lambda_k^{(n)} \in \C^{(m,m)}$ such that
$$
\kern25pt
a(x) P''_n(x) + b(x) P'_n(x) = \sum_{k=-r}^r \Lambda_k^{(n)} P_{n+k}(x)
\qquad \matrix{\hbox{\it(differo-differential} \cr \hbox{ \it equation)}}
$$
\kern25pt with $r\geq \max\{\deg a-2, \deg b-1\}$ independent of $n.$

\medskip

\noindent We use the convention $P_k=0$ for $k<0$.
\et

\bpr See \cite{CMV04,CMV05}. \epr

\br \label{S1R1}
Let us suppose that a $\mathcal{P}$-functional $u\in\PP^{{(m)}'}$ satisfies a
Pearson-type equation $D(u \alpha I) = u\Psi,$ $\alpha\in \mathbb{P}\setminus\{0\},$
$\Psi\in\mathbb{P}^{(m)},$ and let $s = \max\{\deg\alpha-2,\deg\Psi-1\}$. Then, the
proofs given in \cite{CMV05} show that the structure relation appearing in Theorem
\ref{CHSCL} {\it (ii)} is satisfied for the same polynomial $\alpha$ and integer $s.$
However, contrary to the scalar case, the differo-differential equation given in Theorem
\ref{CHSCL} {\it (iii)} can not be ensured for $a=\alpha,$ $r=s,$ but for $a=\alpha^2$
and $r = \max\{2\deg\alpha-2, 2s+2\} = \max\{2\deg\alpha-2, 2\deg\Psi\}\geq s.$
\er

\medskip

In the scalar case, the classical orthogonal polynomials can be characterized by a
Pearson-type equation $D(u\alpha)=u\beta$, $\alpha\in\PP_2\setminus\{0\}$,
$\beta\in\PP_1$, for the corresponding orthogonality functional $u$. When trying to
generalize the concept of classical orthogonal polynomials to the matrix case using a
Pearson-type equation, the following two possibilities appear:
\bi

\item Zero class:
      $u\in \PP^{{(m)}'}$ belongs to the zero class if it is semi-classical with class $s=0$,
      that is, $u$ is quasi-definite and  there exist $\alpha\in\PP_2\setminus\{0\}$,
      $\Psi\in\PP_1^{(m)},$ such that $D(u\alpha I) = u \Psi$.

\item Family $\mathcal{P}_{2,1}$:
      $u \in \PP^{{(m)}'}$ is a $\mathcal{P}_{2,1}$-functional, or belongs to the family
      $\mathcal{P}_{2,1}$, if there exist $\Phi\in\PP_2^{(m)},$  $\Psi\in\PP_1^{(m)},$
      with $\det\Phi\neq0,$ such that $D(u\Phi) = u\Psi.$
\ei

The MOP associated with zero class functionals or quasi-definite
$\mathcal{P}_{2,1}$-functionals can be considered as matrix generalizations of the
classical scalar orthogonal polynomials. Notice that a quasi-definite
$\mathcal{P}_{2,1}$-functional is always semi-classical, but its  class can be greater
than zero. In fact, excepting the scalar case, the family of quasi-definite
$\mathcal{P}_{2,1}$-functionals is strictly greater than the zero class, as can be seen
in Examples 2, 3 and 4. Both, the family $\mathcal{P}_{2,1}$ and the zero class, are
interesting sets of matrix functionals since the related MOP inherit some of the
properties that characterize the classical orthogonal polynomials in the scalar case.
This will be shown in the following sections, which are devoted to the study of the
family $\mathcal{P}_{2,1}$ and the zero class.

\medskip

Before doing that, we will comment some other questions of importance for matrix
orthogonal polynomials. As we have pointed out, a central concept for matrix functionals
is the diagonalizability or, more generally, the reducibility. We say that a functional
$u\in\PP^{{(m)}'}$ is diagonal or block-diagonal if its moment sequence
$(\mu_n)_{n\geq0}$ enjoys such a property. We write $u = u^{(1)} \oplus \cdots \oplus
u^{(k)}$ if $\mu_n = \mu_n^{(1)} \oplus \cdots \oplus \mu_n^{(k)},$ where
$(\mu_n^{(i)})_{n\geq0}$ are the moments of $u^{(i)}.$

To simplify the analysis of a matrix functional $u\in\PP^{{(m)}'}$, the usual strategy is
to connect it with a diagonal or block-diagonal one $\hat u\in\PP^{{(m)}'}$ through a
relation that permits to translate the information from $\hat u$ to $u.$ For instance, if
$\hat u = TuS,$ with $T,S\in\C^{(m,m)}$ non-singular, we say that $u$ is equivalent to
$\hat u.$ In particular, when $S=T^*$ we say that $u$ is congruent to $\hat u$, while if
$S=T^*=T^{-1}$ way say that $u$ is unitarily similar to $\hat u.$ Notice the difference
with the terminology used by other authors, we prefer to preserve the usual one in Linear
Algebra to avoid unnecessary confusion. A matrix functional is diagonalizable or
reducible by equivalence if it is equivalent to a diagonal or block-diagonal one
respectively. We define in a similar way the diagonalizability or reducibility by
congruence and the unitary diagonalizability or reducibility.

A change of variable $t(x)=ax+b,$ $a\in\C\setminus\{0\},$ $b\in\C,$ can be used to relate
matrix functionals too. Given $u\in\PP^{{(m)}'}$ we define $u_t\in\PP^{{(m)}'}$ by
$$
\langle P, u_t \rangle = \langle P \circ t, u \rangle,
$$
so that, if $u=dM$, then $u_t=d(M \circ t^{-1}).$ Notice that, with this definition,
$(Du)_t = (Du_t) t'.$

The kind of relation that we use depends on the properties that we need to preserve. For
example, the equivalence transformation and the change of variable keep invariant the
quasi-definite character, any family $\mathcal{P}_{p,q}$ as well as the class of a
$\mathcal{P}$-functional (in fact, the MOP and the corresponding Pearson-type equations
are trivially related by these transformations). This means that, concerning these
properties, the only non-trivial matrix functionals are those ones that are not reducible
by equivalence or change of variable. In particular, if we are going to study a
characteristic of a functional $u$ that only depends on such properties, then we can
always use the normalization $\langle I, u \rangle = I$ since we can work, for example,
with the equivalent functional $\hat u = u\mu_0^{-1}.$ Also, this allows when studying
zero class functionals to restrict our attention to the canonical choices
$\alpha(x)=1,x,1-x^2,x^2$ of the scalar polynomial in the Pearson-type equation, due to
the freedom in the change of variables.

However, if we are interested in a characteristic that depends on the hermiticity or
positive definiteness of $u$ (or, more generally, on the hermiticity or positive
definiteness of some moments $\mu_n$ or Hankel matrices $\Delta_n$) we must use
congruence transformations and changes of variable with real coefficients. This is the
reason to avoid using the canonical forms of the scalar polynomial $\alpha$ when studying
hermitian zero class functionals, unless we are sure that $\alpha$ has real roots. Also,
the normalization $\langle I, u \rangle = I$ can be used, while preserving any
hermiticity property of $u$, whenever $\mu_0>0$ since, then, we can use the congruent
functional $\hat u = L^{-1}u(L^{-1})^*,$ where $\mu_0=LL^*$ is the Cholesky factorization
of $\mu_0.$

\smallskip

\section{The family $\mathcal{P}_{2,1}$}

\medskip

The aim of this section is to study the differential properties of the MOP associated
with $\mathcal{P}_{2,1}$-functionals. The main result is Theorem \ref{MAIN}, which shows
that some characterizations of the classical scalar orthogonal polynomials remain true
for the matrix family  $\mathcal{P}_{2,1}$. Along the way to prove Theorem \ref{MAIN} we
will obtain a chain of results which have their own interest.

We will start fixing some notations that we will need in the rest of the section. Let
$u\in\PP^{{(m)}'}$ be a $\mathcal{P}_{2,1}$-functional, that is, $D(u\Phi)=u\Psi$, where
$\Phi(x) = \varphi_0 + \varphi_1 x + \varphi_2 x^2,$ $\Psi(x) = \psi_0 + \psi_1 x,$  with
$\varphi_i, \psi_j \in \mathbb{C}^{(m,m)}$ and $\det \Phi\neq 0$. The above Pearson-type
equation is equivalent to
\beq \label{RMM}
n(\mu_{n-1} \varphi_0 + \mu_n\varphi_1 + \mu_{n+1} \varphi_2)
= - (\mu_n \psi_0 + \mu_{n+1} \psi_1),
\qquad n\geq 0,
\eeq
where $(\mu_k)_{k\geq 0}$ are the moments of $u$ and $\mu_{-1} =0$. We denote
$$
\widetilde{u} = u\Phi,\qquad \widetilde{\mu}_n = \langle x^nI, \widetilde{u} \rangle, \qquad
\widetilde{\Delta}_n =
\pmatrix{
\widetilde{\mu}_0&\kern-7pt\widetilde{\mu}_1&\kern-7pt\dots&\kern-7pt\widetilde{\mu}_n
\cr
\dots&\kern-7pt\dots&\kern-7pt\dots&\kern-7pt\dots
\cr
\widetilde{\mu}_n&\kern-7pt\widetilde{\mu}_{n+1}&\kern-7pt&\kern-7pt\widetilde{\mu}_{2n}
}.
$$
The moments of $u$ and $\tilde{u}$ are related by
\beq\label{RMMd}
\widetilde{\mu}_n = \mu_n \varphi_0 + \mu_{n+1} \varphi_1 + \mu_{n+2}\varphi_2,
\quad n\geq 0.
\eeq

One of the characterizations of the classical scalar orthogonal polynomials is that they
are the only sequences of orthogonal polynomials whose derivatives are also sequences of
orthogonal polynomials. The following proposition is the starting point to prove a
similar result for the family $\mathcal{P}_{2,1}$.

\bp\label{FSMOP}
Let $u$ be a $\mathcal{P}_{2,1}$-functional such that $\Delta_0, \Delta_1,\dots,\Delta_n$
are non singular. Then, the corresponding finite segment $(P_k)_{k=0}^n$ of monic MOP
satisfies
$$
\langle x^j P_k',\widetilde{u}\rangle = 0, \qquad j=0,\dots,k-2, \quad k=2,\dots,n,
$$
$$
\langle x^{k-1}P_k', \widetilde{u}\rangle = - E_k(\psi_1 + (k-1)\varphi_2),
\qquad k=1,\dots,n.
$$
\ep

\bpr
From the distributional equation $D(u\Phi) = u\Psi$ we have
$$
\langle x^j P_k, D(u\Phi)\rangle = \langle x^jP_k, u \Psi \rangle,
$$
or, equivalently,
$$
-j\langle x^{j-1}P_k, u\Phi \rangle - \langle x^jP_k', u\Phi \rangle
= \langle  x^jP_k, u\Psi \rangle,
$$
which, for $j=0,\dots, k-1,$ gives the result.
\epr

\bc\label{FSMOPC}
Under the conditions of Proposition \ref{FSMOP}, $(P_k')_{k=1}^n$ is a  finite segment of
MOP with respect to $\widetilde{u}$ if and only if the matrix $\psi_1 + (k-1)\varphi_2$
is non-singular for $k=1,\dots,n.$
\ec

The above corollary shows the interest in finding conditions that ensure the
non-singularity of the matrices $\psi_1 + k \varphi_2$, $k=0,1,2\dots.$ The next lemmas
study the relation between the non-singularity of $\Delta_j$, $j=0,1,\dots,p$, and
$\psi_1 + k \varphi_2$, $k=0,1,\dots,q,$ for small values of $p$ and $q$. They also
inform about the non-singularity of $\tilde{\Delta}_k$, $k=0,1,\dots,q$, a result of
interest since, in the scalar case, $\tilde{u}$ is quasi-definite for any classical
functional $u.$

\bl\label{NNSS}
Let $u$ be a $\mathcal{P}_{2,1}$-functional with $\Delta_0,\Delta_1,\Delta_2$ non-singular.
Then, $\psi_1$ and $\widetilde{\Delta}_0$ are non-singular.
\el

\bpr If $\psi_1$ is singular, there exists $v\in\mathbb{C}^{m}\setminus\{0\}$ such that
$\psi_1 v = 0.$ Relation (\ref{RMM}) for $n=0$ gives $\mu_0\psi_0 + \mu_1 \psi_1 = 0.$
The non-singularity of $\mu_0=\Delta_0$ implies $\psi_0 v=0.$ So, from (\ref{RMM}) we
have
$$
\mu_{n-1}\varphi_0 v +  \mu_n\varphi_1 v +  \mu_{n+1}\varphi_2 v =0,
\qquad  n\geq1,
$$
and, hence,
$$
\Delta_2\pmatrix{\varphi_0v\cr\varphi_1v\cr\varphi_2v\cr}=\pmatrix{0\cr0\cr0\cr}.
$$
Also, $(\varphi_0v,\varphi_1v,\varphi_2v) \neq (0,0,0)$ because $\det\Phi\neq 0$. Now, we
can conclude the singularity of $\Delta_2,$ with contradicts the hypothesis. So, $\psi_1$
is non-singular.

On the other hand, the calculation of $E_1$ gives $E_1=\mu_2-\mu_1\mu_0^{-1}\mu_1,$
which, according to Proposition \ref{S1P1}, is non-singular because $\Delta_1$ is non
singular too. From (\ref{RMMd}) for $n=0$ we get
$$
\widetilde{\mu}_0 = \mu_0 \varphi_0 + \mu_1 \varphi_1 + \mu_2 \varphi_2,
$$
and (\ref{RMM}) for $n=0,1$ gives
$$
\mu_0\psi_0 + \mu_1\psi_1 = 0, \qquad
\mu_0\varphi_0+\mu_1\varphi_1+\mu_2\varphi_2 = -(\mu_1\psi_0+\mu_2\psi_1).
 $$
Therefore,
$$
\widetilde{\Delta}_0 = \widetilde{\mu}_0 = - \mu_1 \psi_0 - \mu_2 \psi_1 =
-(\mu_2-\mu_1\mu_0^{-1}\mu_1)\psi_1=-E_1\psi_1
$$
is non-singular.
\epr

As a first consequence, we obtain the following announced result.

\bt\label{MC}
If $u\in \mathcal{P}_{2,1}$ and $\Delta_0, \Delta_1, \Delta_2$ are non-singular,
the $\C^{(m,m)}$-right-module $\mathcal{M}_{2,1}(u)$ is cyclic.
\et

\bpr
Let us suppose that $D(u \Phi^{(i)}) = u \Psi^{(i)}$ with $\Phi^{(i)}\in\PP_2^{(m)},$
$\Psi^{(i)}\in \PP_1^{(m)}$ for $i=1,2,$ and assume that $\det \Phi^{(1)} \neq 0.$ We are
going to prove that $\Phi^{(2)} = \Phi^{(1)}\Lambda,$  $\Lambda \in \mathbb{C}^{(m,m)}.$
Let $\Psi^{(i)}(x)= \psi_0^{(i)} + \psi_1^{(i)} x$ with
$\psi_0^{(i)},\psi_1^{(i)}\in\mathbb{C}^{(m,m)}.$ Since $u$ satisfies the hypothesis of
Lemma \ref{NNSS}, $\psi_1^{(1)}$ is non-singular. Hence,
$A=\Phi^{(1)}(\psi_1^{(1)})^{-1}\psi_1^{(2)}-\Phi^{(2)}$ satisfies
$$
D(uA) = u \left( \psi_0^{(1)} (\psi_1^{(1)})^{-1} \psi_1^{(2)} - \psi_0^{(2)} \right).
$$
From (\ref{RMM}) for $n=0$, $\psi_0^{(i)} = -\mu_0^{-1}\mu_1\psi_1^{(i)},$ therefore,
$D(uA)=0$. If $A(x) = A_0+A_1x+A_2x^2,$ $A_i\in \mathbb{C}^{(m,m)},$ we get $\mu_n A_0 +
\mu_{n+1}A_1 + \mu_{n+2}A_2=0$ for $n\geq 0,$ which implies
$$
\Delta_2 \pmatrix{A_0\cr A_1\cr A_2\cr}= 0.
$$
Since $\Delta_2$ is non-singular, $A=0$ and, thus,
$\Phi^{(2)} = \Phi^{(1)}\bigl(\psi_1^{(1)}\bigr)^{-1} \psi_1^{(2)}.$
\epr

Now, we are going to consider  $\mathcal{P}_{2,1}$-functionals satisfying the hypothesis
of Lemma \ref{NNSS}. In such a case we can write $\psi_1= I$ without loss of generality
because the Pearson-type equation can be written as
$D(u\Phi\psi_1^{-1})=u\Psi\psi_1^{-1}$.

\bl\label{NNSSS}
Let $u$ be a $\mathcal{P}_{2,1}$-functional with $\Delta_k$ non-singular
for $k=0,1,2,3.$ Then,

\item{(i)} $\psi_1$ and $\psi_1 + \varphi_2$ are non-singular.

\item{(ii)} $\widetilde{\Delta}_0$ and $\widetilde{\Delta}_1$ are non-singular.

\item{(iii)} $\tilde{u}$ is a $\mathcal{P}_{2,1}$-functional, that is,
$D(\widetilde{u}\widetilde{\Phi})=\widetilde{u}\widetilde{\Psi}$, with
$\widetilde{\Phi}(x) =\sum_{i=0}^2\tilde{\varphi_i} x^i$,
${}\kern22pt \widetilde{\Psi}(x) = \sum_{j=0}^1 \tilde{\psi}_j x^j$, where
$\tilde{\varphi}_i, \tilde{\psi}_j \in \mathbb{C}^{(m,m)}$ and
$\det\tilde{\Phi}\neq 0$.
Moreover, ${}\kern22pt \tilde{\Phi}$, $\tilde{\Psi}$ can be chosen such that
$\tilde{\varphi_2}=\psi_1^{-1}\varphi_2$ and $\tilde{\psi}_1=\psi_1^{-1}(\psi_1+2\varphi_2)$.
\el

\bpr
We will assume without of loss of generality that $\psi_1 = I.$

\item{\it (i)} Let us suppose that $I+\varphi_2$ is singular. There exists
$v\in\mathbb{C}^{m}\setminus\{0\}$ such that $\varphi_2 v = -v.$ Writing (\ref{RMM})
for $n=0,1,$
$$
\mu_1 +\mu_0\psi_0 = 0, \qquad \mu_1(\psi_0+\varphi_1)v + \mu_0 \varphi_0 v = 0.
$$
Then,
\beq\label{RMM0}
-\psi_0(\psi_0+\varphi_1)v + \varphi_0 v = 0.
\eeq

Consider (\ref{RMM}) again, but for $n$ and $n+1$:
$$
\cases{
n\mu_{n-1}\varphi_0 + \mu_n(\psi_0 + n\varphi_1) + \mu_{n+1} (I + n \varphi_2) = 0,
\smallskip
\cr
(n+1) \mu_n\varphi_0 + \mu_{n+1}[\psi_0 + (n+1) \varphi_1] + \mu_{n+2} [I + (n+1)\varphi_2] = 0.
}
$$
Multiplying the first equation on the right by $\psi_0 + \varphi_1$ and subtracting
the second one, gives
$$
\kern-20pt
n \mu_{n-1}\varphi_0(\psi_0 + \varphi_1)
+ \mu_n\left[
\psi_0\left(\psi_0 + \varphi_1\right)
- \varphi_0 + n  \left(\varphi_1\left(\psi_0 + \varphi_1\right)
- \varphi_0\right)
\right] +
$$
$$
\kern93pt
+ \; n\mu_{n+1} \left[\varphi_2(\psi_0 + \varphi_1) - \varphi_1\right]
- \mu_{n+2}\left[I+(n+1)\varphi_2\right] = 0.
$$
Then, taking into account (\ref{RMM0}), we get
$$
\kern-120pt
\mu_{n-1}\varphi_0\left(\psi_0 + \varphi_1\right) v
+ \mu_n \left[\varphi_1\left(\psi_0 + \varphi_1\right) -\varphi_0\right] v  +
$$
\beq \label{RMMM}
\kern77pt
+ \; \mu_{n+1}\left[\varphi_2\left(\psi_0 + \varphi_1\right)
- \varphi_1\right] v - \mu_{n+2} v = 0, \quad n\geq 1,
\eeq
which implies
$$
\Delta_3
\pmatrix{
\varphi_0(\psi_0 + \varphi_1)v
\cr
[\varphi_1(\psi_0 + \varphi_1) - \varphi_0]v
\cr
[\varphi_2(\psi_0 + \varphi_1) - \varphi_1]v
\cr
-v}
= \pmatrix{0\cr0\cr0\cr0}.
$$
This contradicts the non-singularity of $\Delta_3.$

\item{\it (ii)} By Proposition \ref{FSMOP} and Corollary \ref{FSMOPC}, $\{P_1', P_2'\}$
is a finite segment of MOP with respect to $\widetilde{u}.$ The result follows from
Proposition \ref{S1P1}.

\item{\it (iii)} The existence of matrix polynomials $\widetilde{\Phi}, \widetilde{\Psi}$
satisfying $D(\tilde{u}\tilde{\Phi})= \tilde{u}\tilde{\Psi}$ is ensured if
\beq \label{EDM}
\Psi \widetilde{\Phi} + \Phi\widetilde{\Phi}' = \Phi\widetilde{\Psi}.
\eeq
Writing
$
\widetilde{\Phi}(x) =
\widetilde{\varphi}_0 + \widetilde{\varphi}_1 x + \widetilde{\varphi}_2 x^2,
$
$
\widetilde{\Psi}(x) = \widetilde{\psi}_0 + \widetilde{\psi}_1 x,
$
(\ref{EDM}) is equivalent to the system
\beq\label{EDMM}
   \pmatrix{\psi_0&0&\varphi_0&0\cr
            I&\psi_0&\varphi_1&0\cr
            0&I&\varphi_2&0\cr
            0&0&0&I+2\varphi_2}
   \pmatrix{\widetilde{\varphi}_0\cr \widetilde{\varphi}_1\cr
            \widetilde{\varphi}_1-\widetilde{\psi}_0\cr\widetilde{\varphi}_2} =
   \pmatrix{0\cr\varphi_0(\widetilde{\psi}_1-2\widetilde{\varphi}_2)\cr
            \varphi_1(\widetilde{\psi}_1-2\widetilde{\varphi}_2)-\psi_0\widetilde{\varphi}_2\cr
            \varphi_2\widetilde{\psi}_1\cr}.
\eeq
A solution of the last equation is
$\widetilde{\psi}_1 = I + 2\varphi_2$,
$\widetilde{\varphi}_2 = \varphi_2$.
With this choice, converting the system into triangular form gives
$$
\kern-3pt
     \pmatrix{I&\psi_0&\varphi_1\cr
              0&I&\varphi_2\cr
              0&0&\varphi_0-\psi_0\varphi_1+\psi_0^2\varphi_2}
\kern-3pt
     \pmatrix{\widetilde{\varphi}_0\cr
              \widetilde{\varphi}_1\cr
              \widetilde{\varphi}_1-\widetilde{\psi}_0}
\kern-2pt = \kern-1pt
     \pmatrix{\varphi_0\cr
              \varphi_1-\psi_0\varphi_2\cr
              -\psi_0(\varphi_0-\psi_0\varphi_1+\psi_0^2\varphi_2)}.
$$
From (\ref{RMM}) for $n=0,$ $\mu_0 \psi_0 + \mu_1 = 0$, so,
$$
\Upsilon := \varphi_0 - \psi_0\varphi_1 + \psi_0^2 \varphi_2
= \varphi_0 + \mu_0^{-1}\mu_1\varphi_1 + (\mu_0^{-1}\mu_1)^2\varphi_2 =
$$
$$
= \mu_0^{-1}(\mu_0\varphi_0  + \mu_1\varphi_1 + \mu_1\mu_0^{-1}\mu_1\varphi_2).
$$
Since $E_1=\mu_2 - \mu_1\mu_0^{-1}\mu_1$,
$$
\Upsilon =
\mu_0^{-1}(\mu_0\varphi_0 + \mu_1\varphi_1 + \mu_2\varphi_2 - E_1\varphi_2)
$$
that, keeping in mind (\ref{RMM}) for $n=1$, can be expressed as
$$
\kern-100pt
\Upsilon = - \mu_0^{-1}(\mu_1\psi_0 + \mu_2 + E_1\varphi_2) =
$$
$$
\kern22pt
= - \mu_0^{-1}(-\mu_1\mu_0^{-1}\mu_1 + \mu_2 + E_1\varphi_2)
= - \mu_0^{-1} E_1(I + \varphi_2).
$$
That is, $\Upsilon$ is non-singular, what ensures that (\ref{EDMM}) has a solution.

Finally, we are going to prove that $\det\widetilde{\Phi} \neq 0.$
From (\ref{EDM}) we can deduce
$$
\Phi(\tilde{\Psi} -\tilde{\Phi}') = \Psi\tilde{\Phi}.
$$
Since $\det\Phi\neq0,$ $\det\tilde\Phi=0$ implies $\det(\tilde{\Psi}-\tilde{\Phi}')=0$.
However, taking into account that $\tilde{\psi}_1 = I + 2\varphi_2$ and
$\tilde{\varphi}_2 = \varphi_2$ we get $\tilde{\Psi}(x)-\tilde{\Phi}'(x) = \tilde{\psi}_0
- \tilde{\varphi}_1 + I x,$ which has non-null determinant.
\epr

\bl\label{NNSSSS}
Let $u$ be a $\mathcal{P}_{2,1}$-functional with $\Delta_k$ non-singular for $k=0,1,2,3,4.$
Then,

\item{(i)}  $\psi_1+j \varphi_2$ is non-singular for $j=0,1,2$.

\item{(ii)} $\widetilde{\Delta}_j$ is non-singular  for $j=0,1,2$.
\el

\bpr
We will assume without of loss of generality that $\psi_1 = I.$

\item{\it (i)} Taking into account Lemma \ref{NNSSS} {\it (iii)}, the functional $\widetilde{u}$
satisfies $D(\widetilde{u}\widetilde{\Phi})=\widetilde{u}\widetilde{\Psi},$ with
$\widetilde{\varphi}_2 = \varphi_2,$ $\;\widetilde{\psi}_1 = I + 2\varphi_2,$ where
$\;\tilde{\varphi}_i,\ \tilde{\psi}_j$ have the same meaning as in the proof of the previous lemma.

Let us suppose that $I + 2\varphi_2$ is singular. Then, there exists
$v\in\mathbb{C}^{m}\setminus\{0\}$ such that $\varphi_2 v = -{1\over 2} v,$ that is,
$\widetilde{\psi}_1 v = 0.$ Since
$D(\widetilde{u}\widetilde{\Phi})=\widetilde{u}\widetilde{\Psi},$ we have
$$
n(\widetilde{\mu}_{n-1}\widetilde{\varphi}_0
+ \widetilde{\mu}_n\widetilde{\varphi}_1
+ \widetilde{\mu}_{n+1}\widetilde{\varphi}_2) =
-(\widetilde{\mu}_n \widetilde{\psi}_0
+ \widetilde{\mu}_{n+1}\widetilde{\psi}_1),
\qquad n\geq 0,
$$
which, for $n=0,$ gives
$\widetilde{\mu}_0\widetilde{\psi}_0+\widetilde{\mu}_1\widetilde{\psi}_1 = 0.$  Hence,
$\widetilde{\psi}_0 v= 0$ because, from Lemma \ref{NNSS},
$\widetilde{\mu}_0=\tilde{\Delta}_0$ is non-singular. So,
\beq \label{MUTILDE}
(\widetilde{\mu}_{n-1}\widetilde{\varphi}_0
+ \widetilde{\mu}_n\widetilde{ \varphi}_1
+ \widetilde{\mu}_{n+1}\widetilde{\varphi}_2) v = 0,
\quad n\geq 1.
\eeq
According to
(\ref{RMMd}),
 $$
 \kern-35pt
 \mu_{n-1}\varphi_0\widetilde{\varphi}_0 v
 + \mu_n(\varphi_1 \widetilde{\varphi}_0
 + \varphi_0 \tilde{\varphi_1}) v
 + \mu_{n+1}(\varphi_2\widetilde{\varphi}_0
 + \varphi_1\widetilde{\varphi}_1
 + \varphi_0\widetilde{\varphi}_2)v +
 $$
 $$
 \kern114pt
 + \;\mu_{n+2}(\varphi_2\widetilde{\varphi}_1
 + \varphi_1\widetilde{\varphi}_2)v
 + \mu_{n+3}\varphi_2\widetilde{\varphi}_2 v = 0,\quad n\geq 1,
 $$
and from here we can deduce the singularity of $\Delta_4,$ because
$\varphi_2 \widetilde{\varphi}_2 v = \varphi_2^2 v = {1\over 4} v \neq 0.$
This contradicts the hypothesis. So, $\widetilde{\psi}_1$ is non-singular.

\item{\it (ii)} From Corollary \ref{FSMOPC}, $\{P_1', P_2', P_3'\}$ is a finite
segment of MOP with respect to $\widetilde{u}$ and, so, Proposition \ref{S1P1}
ensures that $\widetilde{\Delta}_2$ is non-singular.
\epr

The previous lemmas can be generalized through an inductive process. This process will
need the following result too.

\bl\label{GENERAL}
Let $u\in\mathbb{P}^{(m)'}$ and $F\in\mathbb{P}_p^{(m)}$, with $\det F\neq 0.$ We denote
$\widetilde{u} = u F$ and we suppose that there exist $v_0,v_1,\dots,v_q\in\mathbb{C}^m$,
with $v_k\neq 0$ for some $k \in \{0,1,\dots, q\},$ such that the moments
$(\widetilde{\mu}_n)_{n\geq 0}$  of the functional $\widetilde{u}$ satisfy
$$
\sum_{j=0}^q \widetilde{\mu}_{n+j} v_j = 0,\qquad \forall n\geq 0.
$$
Then, there exist $ w_1, w_2,\dots, w_{p+q}\in \mathbb{C}^{m},$ with $\;w_k\neq 0$ for
some $\;k\in \{0,\dots, p+q\}$, such that  the moments $(\mu_n)_{n\geq 0}$ of the
functional $u$ satisfy
$$
\sum_{k=0}^{p+q} \mu_{n+k} w_k = 0, \qquad \forall n\geq 0.
$$
\el

\bpr
We will write $F(x) = f_0 + f_1  x + \cdots + f_p x^p$ with $f_i\in\mathbb{C}^{(m,m)}.$ Then,
$\widetilde{\mu}_n = \sum_{i=0}^p \mu_{n+i}f_i$ and the hypothesis of the lemma gives
$$
0 = \sum_{j=0}^q \widetilde{\mu}_{n+j} v_j =
    \sum_{j=0}^q \left(\sum_{i=0}^p \mu_{n+j+i} f_i\right)v_j =
    \sum_{k=0}^{p+q} \mu_{n+k}\sum_{i=0}^p f_i v_{k-i},
$$
with the convention $v_{-1} = \cdots =v_{-p}=0.$ So, the vectors
$w_k = \sum_{i=0}^p f_i v_{k-i},$ $k=0,\dots,p+q,$ satisfy the equality of the statement.
It will be enough to  prove that not all the vectors $w_k$ are null. If all of them are zero,
$\sum_{i=0}^p f_i v_{k-i} = 0$ for $k=0, \dots,p+q,$ and this implies
$$
0  = \sum_{k=0}^{p+q} x^k \sum_{i=0}^p f_i v_{k-i}, \qquad \forall x \in \mathbb{C},
$$
or, equivalently,
$$
0 = \sum_{j=0}^q x^j\left(\sum_{i=0}^p f_i x^i\right) v_j =
    F(x) \sum_{j=0}^q v_j x^j, \qquad \forall x \in \mathbb{C}.
$$
Since $\det F\neq 0,$ we obtain from Remark \ref{DET} that $\sum_{j=0}^q v_j x^j = 0$
for all $x\in\C$, which means that $v_j=0$ for $j=0,\dots,q,$ in contradiction with the hypothesis.
\epr

Now we can reach the generalization of Lemmas \ref{NNSS}, \ref{NNSSS} and \ref{NNSSSS}.

\bt\label{SFS}
Let $u$ be  a  $\mathcal{P}_{2,1}$-functional with $\Delta_k$ non-singular for
$k=0,1,\dots,n,$ where $n\geq 2.$  Then, $\psi_1 + j\varphi_2 $ and $\tilde{\Delta}_j$
are non-singular for $j=0,1,\dots,n-2$.
\et

\bpr
Due to Lemmas \ref{NNSS}, \ref{NNSSS} and \ref{NNSSSS} the result is true for
$n=2, 3, 4.$ We will assume the statement for an index $n\geq 2$, and we will prove that
it is also true for $n+1$.

Assume that $\Delta_0, \Delta_1,\dots,\Delta_n, \Delta_{n+1}$ are non-singular. Then, the
hypothesis of induction implies that $\psi_1 + j \varphi_2$ and $\tilde{\Delta}_j$ are
non-singular for $j=0,1,\dots,n-2$. We only must prove that $\psi_1+(n-1)\varphi_2$ and
$\tilde{\Delta}_{n-1}$ are non-singular too. For this purpose we will introduce a set of
$\mathcal{P}_{2,1}$-functionals $u^{(j)}$, $j=0,1,\dots,$  using the superscript $(j)$
for the associated elements.

Let us define $u^{(0)} = u,$ $\Phi^{(0)} = \Phi,$ $\Psi^{(0)}=\Psi$. Taking into account
Lemmas \ref{NNSSS} and \ref{NNSSSS}, given
$u^{(1)}=u^{(0)}\Phi^{(0)}\bigl(\psi_1^{(0)}\bigr)^{-1}$ there exist
$\Phi^{(1)}\in\mathbb{P}_2^{(m)},$  $\Psi^{(1)}\in\mathbb{P}_1^{(m)},$ satisfying
$D\left(u^{(1)} \Phi^{(1)}\right) = u^{(1)} \Psi^{(1)},$ with $\det \Phi^{(1)}\neq 0,$
$\varphi_2^{(1)}=\varphi_2^{(0)}$ and $\psi_1^{(1)}=\psi_1^{(0)}+2\varphi_2^{(0)}$
non-singular. Moreover, from Proposition \ref{FSMOP}, $E_k^{(1)} = -{1\over k+1}
E_{k+1}^{(0)} \bigl(\psi_1^{(0)} + k \varphi_2^{(0)}\bigr).$ This implies that
$E_0^{(1)},\dots,E_{n-2}^{(1)}$ and, thus, $\Delta_0^{(1)},\dots,\Delta_{n-2}^{(1)}$ are
non-singular.

Following this procedure, we can construct inductively a set of
$\mathcal{P}_{2,1}$-functionals $u^{(j)}$, $j=0,1,\dots,l-1$
$\left(l=\left[{n\over2}\right]\right),$ satisfying
$$
u^{(j+1)} = u^{(j)} \Phi^{(j)} \bigl(\psi_1^{(j)}\bigr)^{-1},
$$
$$
D(u^{(j)} \Phi^{(j)}) = u^{(j)} \Psi^{(j)}, \quad \varphi_2^{(j)}=\varphi_2, \quad
\psi_1^{(j)} = \psi_1 + 2 j \varphi_2,
$$
$$
E_k^{(j+1)} = -{1\over k+1} E_{k+1}^{(j)} \left[\psi_1 + (2j+k)\varphi_2\right],
$$
$$
\Delta_0^{(j)},\dots,\Delta_{n-2j}^{(j)} \quad \hbox{non-singular}.
$$

Let us suppose that $n$ is even $(n=2l)$. Then,
$\Delta_0^{(l-1)},\Delta_1^{(l-1)},\Delta_2^{(l-1)}$ are non-singular. If
$\psi_1+(n-1)\varphi_2=\psi_1^{(l-1)} + \varphi_2^{(l-1)}$ is singular, the same
arguments that lead to (4) in the proof of Lemma \ref{NNSSS} give now
$$
\sum_{j=0}^3 \mu_{k+j}^{(l-1)} v_j=0, \qquad v_3\neq 0, \qquad k\geq 0.
$$
Since $u^{(l-1)} = u F$, $\deg F \leq 2l-2 = n-2$, we get from Lemma \ref{GENERAL}
$$
\sum_{j=0}^{n+1} \mu_{k+j} w_j = 0, \qquad  \hbox{some}\;w_j\neq 0, \qquad k\geq 0.
$$
This contradicts the non-singularity of $\Delta_{n+1}$, so, $\psi_1 + (n-1)\varphi_2$
must be non-singular.

If, on the contrary, $n$ is odd $(n=2l+1)$,
$\Delta_0^{(l-1)},\Delta_1^{(l-1)},\Delta_2^{(l-1)},\Delta_3^{(l-1)}$ are non-singular.
Thus, analogously to (8) in the proof of Lemma \ref{NNSSSS}, we find that, if
$\psi_1+(n-1)\varphi_2=\psi_1^{(l-1)} + 2\varphi_2^{(l-1)}$ is singular,
$$
\sum_{j=0}^4 \mu_{k+j}^{(l-1)} v_j = 0, \qquad v_4\neq 0, \qquad k\geq 0.
$$
Now, $u^{(l-1)} = u F$, $\deg F \leq 2l-2 = n-3$, so, Lemma \ref{GENERAL} gives again
the same condition
$$
\sum_{j=0}^{n+1} \mu_{k+j} w_j = 0, \qquad \hbox{some}\;w_j\neq 0, \qquad k\geq 0,
$$
so, $\psi_1+(n-1)\varphi_2$ is also non-singular in this case.

Finally, the non-singularity of $\tilde{\Delta}_{n-1}$ follows from Proposition
\ref{S1P1} and the relation $\tilde{E}_{n-1} = -{1\over n} E_n\bigl(\psi_1 +
(n-1)\varphi_2\bigr)$ given in Proposition \ref{FSMOP}.
\epr

The previous theorem and Corollary \ref{FSMOPC} have the following immediate consequences.

\bc\label{SSC1}
If $u$ is a quasi-definite $\mathcal{P}_{2,1}$-functional, then $\psi_1+n\varphi_2$
is non-singular for $n=0,1,2,\dots$.
\ec

\bc\label{SSC2}
If $u$ is a quasi-definite $\mathcal{P}_{2,1}$-functional, then $\widetilde{u} = u \Phi$
is a quasi-definite $\mathcal{P}_{2,1}$-functional too. Moreover, if $(P_n)_{n\geq0}$ is
the sequence of monic MOP with respect to $u$, then $\left({1\over n}
P_n'\right)_{n\geq1}$ is the sequence of monic MOP with respect to $\tilde{u}.$
\ec

\br\label{SSR}
The Pearson-type equation $D(u \Phi)=u\Psi,$ $\Phi\in\mathbb{P}_2^{(m)},$
$\Psi\in\mathbb{P}_1^{(m)},$ is equivalent to the recurrence
$n\mu_{n-1}\varphi_0+\mu_n(\psi_0+n\varphi_1)+\mu_{n+1}(\psi_1+n\varphi_2)=0,$ $n\geq0$.
Therefore, the non-singularity of the matrices $\psi_1+n\varphi_2$ for $n\geq0$ is a
sufficient condition for the existence of a solution $u$ of the Pearson-type equation.
Indeed, this condition ensures that the solutions are determined by $\mu_0 = \langle I,u
\rangle$ or, in other words, the solution is unique up to left matrix factors. Then,
according to Corollary \ref{SSC1}, if the Pearson-type equation has a quasi-definite
solution, the quasi-definite solutions are exactly those solutions determined by a
non-singular matrix $\mu_0$.
\er

\subsection{Characterization of  the  family $\mathcal{P}_{2,1}$}

\smallskip

In the scalar case, the classical orthogonal polynomials can be characterized
alternatively by a Pearson-type equation (see \cite{Ch78,Mar91,Mar93,Sho39}), the
orthogonality of the derivatives (see \cite{BLN87,Ch78,HVR85,Mar91,Mar93}) or a linear
relation between the polynomials $P_n$ and $P_{n+1}',$ $P_n',$ $P_{n-1}'$ (see
\cite{MBP94}). The consequences of the previous analysis provide an analogue of these
equivalences for the matrix case, which constitute a characterization of the
quasi-definite $\mathcal{P}_{2,1}$-functionals. In the proof of this characterization we
will need the following results too.

\bl\label{NULL}
Let $u \in \mathbb{P}^{(m)'}$ such that $\Delta_n$ is non-singular. Then,
$$
uP=0, \;\, P\in\PP_n^{(m)} \;\;\Rightarrow\;\; P=0.
$$
\el

\bpr
Let $P(x)=\sum_{i=0}^nA_ix^i$, $A_i\in\C^{(m,m)}.$ Then, $uP=0$ is equivalent to
$\mu_kA_0+\cdots+\mu_{k+n}A_n=0$ for $k\geq0$, which implies
$$
\Delta_n \pmatrix{A_0 \cr \vdots \cr A_n} = 0,
$$
and, thus, $P=0$ if $\Delta_n$ is non-singular.
\epr

\bp\label{QOMOP}
Let $u, v \in \mathbb{P}^{(m)'}$ with $u$ quasi-definite and $(P_n)$ its corresponding
sequence of monic MOP. Then, the following statements are equivalent:

    \item{(i)} $v=uA, \ \ A\in \mathbb{P}_p^{(m)}.$

    \item{(ii)} $(P_n)$ is quasi-orthogonal of order not greater than $p$
                with respect to $v$: \vskip3pt
                $\kern3pt \langle x^k P_n, v \rangle = 0,  \ \ k=0,\dots,n-p-1.$
\ep

\bpr See \cite{CMV05}. \epr

Here is the referred characterization of the quasi-definite $\mathcal{P}_{2,1}$-functionals.

\bt\label{MAIN}
Let $u \in \mathbb{P}^{(m)'}$ be quasi-definite and $(P_n)$ its sequence of monic MOP.
Then, the following assertions are equivalent:

\item {(i)} $u$ is a  $\mathcal{P}_{2,1}$-functional.

\item{(ii)} $\bigl(P_n'\bigr)$ is a sequence of  MOP with respect
to a quasi-definite  functional $\widetilde{u}.$

\item{(iii)} There exist matrices $a_n, b_n \in \mathbb{C}^{(m,m)}$ such that
$$
\kern25pt
P_n={1\over n+1} P_{n+1}' + a_n P_n' + b_n P_{n-1}', \qquad n\geq 0,
$$
\kern25pt with $\gamma_n - b_n$ non-singular for $n \geq 1.$

\noindent Moreover, $\widetilde{u} = u \Phi$, $\Phi \in \mathbb{P}_2^{(m)},$ $\det\Phi\neq0$
and $D(u \Phi) = u \Psi,$ $\Psi\in \mathbb{P}_1^{(m)}.$
Besides, $\widetilde{u}$ is a quasi-definite $\mathcal{P}_{2,1}$-functional too.
\et

\bpr
\item{\it (ii) $\Leftrightarrow$ (iii)}
The sequence of matrix polynomials $(P_n)$ satisfies the recurrence relation,
$$
x P_n = P_{n+1} + \beta_n P_n + \gamma_n P_{n-1},
$$
so,
\beq\label{TH1}
P_n = - x P_n' + P_{n+1}' + \beta_n P_n' + \gamma_n P_{n-1}'.
\eeq
If we assume {\it (ii)}, $(P_n')$ also satisfies a recurrence relation
\beq\label{TH2}
{1\over n} x P_n' =
  {1\over n+1} P_{n+1}'
+ {1\over n} \tilde{\beta}_{n-1} P_n'
+ {1\over n-1} \tilde{\gamma}_{n-1} P_{n-1}'
\eeq
and, then, (\ref{TH1}) and (\ref{TH2}) imply
\beq\label{TH3}
P_n = {1\over n+1} P_{n+1}'+ a_n P_n' + b_n P_{n-1}',
\eeq
where $a_n=\beta_n-\tilde{\beta}_{n-1}$ and $b_n=\gamma_n-{n\over n-1}\tilde{\gamma}_{n-1}$.
Notice that $\gamma_n-b_n={n\over n-1}\tilde{\gamma}_{n-1}$ is non-singular.

For the converse, from (\ref{TH1}) and (\ref{TH3}),
$$
{1\over n} x P_n' =
  {1\over n+1} P_{n+1}'
+ {1\over n}\left(\beta_n - a_n\right) P_n'
+ {1\over n} \left(\gamma_n - b_n\right)P_{n-1}'.
$$
Now, we have a recurrence relation for $(P_n')$ with $\tilde{\beta}_{n-1}=\beta_n-a_n$
and $\tilde{\gamma}_{n-1}={n-1\over n}(\gamma_n-b_n).$ Since $\gamma_n-b_n$ is
non-singular, the Favard theorem assures the existence of a functional
$\tilde{u}\in\mathbb{P}^{(m)'}$ such that $(P_n')$ is a sequence of  MOP with respect to
$\tilde{u}.$

\smallskip

\item{\it (ii), (iii) $\Rightarrow$ (i)}
Assume the  relation $P_n = {1\over n+1}P_{n+1}'+a_n P_n'+b_n P_{n-1}'$ and the fact that
$(P_n')$ is a sequence of MOP with respect to a certain functional $\tilde{u}.$ Notice
that this last hypothesis implies the non-singularity of $\tilde{E}_{n-1} = {1\over
n}\langle x^{n-1} P_n', \tilde{u}\rangle$ for $n\geq1$. Under the assumptions,
$$
\langle x^k P_n, \tilde{u}\rangle = 0, \qquad k=0, \dots, n-3.
$$
So, $(P_n)$ is a quasi-orthogonal sequence with respect to $\tilde{u}$ of order not
greater than 2. Proposition \ref{QOMOP} says that there exists
$\Phi\in\mathbb{P}_2^{(m)}$ such that $\tilde{u}=u\Phi.$ Setting $w = D(u\Phi)$,
$$
\langle x^k P_n, w \rangle =
- k \langle x^{k-1}P_n, u\Phi\rangle -\langle x^k P_n', u\Phi\rangle = 0,
\quad k=0,\dots,n-2.
 $$
Hence, $(P_n)$ is quasi-orthogonal with respect to $w$ of order not greater than 1 and,
thus, there exists $\Psi\in \mathbb{P}_1^{(m)}$ such that $w = u\Psi.$

It only remains to prove that $\det\Phi\neq0$. For this purpose, notice that the equality
$$
\langle x^{n-1}P_n, D(u\Phi)\rangle = \langle x^{n-1}P_n, u\Psi\rangle
$$
gives
$$
-(n-1) E_n \varphi_2 - n \tilde{E}_{n-1} = E_n \psi_1.
$$
Hence, $\psi_1 + (n-1) \varphi_2$ is non-singular for $n\geq 1$. Suppose $\det\Phi=0$.
Then, according to Remark \ref{DET}, there exists $v\in\C^{m}[x]\setminus\{0\}$ such that
$\Phi v = 0$. Consider the matrix polynomial $A\in\PP^{(m)}$ whose columns are all equal
to $v$. Taking into account Lemma \ref{NULL}, the equality
$$
u (\Psi-\Phi') A = (Du)\Phi A = 0
$$
proves that $(\Psi-\Phi')v=0$. So, $\Psi v + \Phi v' = 0$ and, if
$v(x)=v_0+\cdots+v_nx^n$, $v_i\in\C^{m}$, with $v_n\neq0$, we get
$(\psi_1+n\varphi_2)v_n=0$, which is impossible.

\smallskip

\item{\it (i) $\Rightarrow$ (ii)} This implication is given by Corollary \ref{SSC2}.
\epr

\br \label{CHAIN}
Theorem \ref{MAIN} ensures that any quasi-definite $\mathcal{P}_{2,1}$-functional $u$
generates a sequence $(u^{(n)})_{n\geq 0}$ of quasi-definite
$\mathcal{P}_{2,1}$-functionals, starting with $u^{(0)}=u$, and such that, for
$n\geq0$,
$$
u^{(n+1)}=u^{(n)}\Phi^{(n)},
\quad \Phi^{(n)}\in\mathbb{P}_2^{(m)},
\quad \det\Phi^{(n)}\neq0,
$$
$$
\quad D(u^{(n)}\Phi^{(n)})= u^{(n)}\Psi^{(n)},
\quad \Psi^{(n)}\in\mathbb{P}_1^{(m)}.
$$
Moreover, the $k$-th derivatives $(P_n^{(k)})_{n\geq k}$  form a sequence of MOP with
respect to $u^{(k)}.$ That is, as in the scalar case, if the first derivatives of a
sequence of MOP are orthogonal, the higher order derivatives are orthogonal too.
\er

\vskip0pt

\br
If $u$ is not quasi-definite but $\Delta_0,\dots,\Delta_n$ are non-singular, {\it (ii)}
and {\it (iii)} remain equivalent, but for the finite segment $\left(P_k\right)_{k=0}^n$
of monic MOP with respect to $u.$ Besides, in this case, {\it (i)} also implies {\it
(ii)} and {\it (iii)}, but only for the finite segment $\left(P_k\right)_{k=0}^{n-1},$
according to Theorem  \ref{SFS}.
\er

\medskip

The following consequence of Theorem \ref{MAIN} will be of interest when studying the
differential equation associated with the zero class MOP.

\bc\label{RRPPS}
If a sequence $\left(P_n\right)$ of monic MOP belongs to the family $\mathcal{P}_{2,1}$,
then $P_{n\pm 1}' \in {\rm span}_{\mathbb{C}^{(m,m)}}\{x P_n', P_n', P_n\}$. More
precisely,
$$
\kern-65pt
P_{n-1}' = E_{n-1} M_{n-2} M_{2n-1}^{-1}E_n^{-1}\left\{\left(x+{1\over n}\pi_n\right)P_n'
- n P_n\right\},
$$
$$
\kern-30pt
P_{n+1}' = (n+1) E_n \left\{\left(\varphi_2 M_{2n-1}^{-1}E_n^{-1} x - {1\over n}
M_{2n-2}M_{2n-1}^{-1} E_n^{-1} \pi_n + \right.\right.
$$
$$
\kern100pt
\left.\left. +{1\over n+1}
E_n^{-1} \pi_{n+1}\right) P_n' + M_{n-1}M_{2n-1}^{-1} E_n^{-1} P_n\right\},
$$
where $E_n = \langle x^n P_n, u\rangle,$ $P_n(x)=x^n+\pi_nx^{n-1}+\cdots$ and
$M_n=\psi_1+n\varphi_2.$
\ec

\bpr

Using (\ref{TH1}) and (\ref{TH3}), we get by eliminating $P_{n+1}'$ and  $P_{n-1}'$
respectively,
$$
\cases{
nP_n = \left(x-\beta_n+(n+1)a_n\right)P_n' - \left(\gamma_n-(n+1)b_n\right)P_{n-1}',
\smallskip
\cr
\left(1-b_n\gamma_n^{-1}\right)P_n =
\left({1\over n+1}-b_n\gamma_n^{-1}\right)P_{n+1}' +
\left(b_n\gamma_n^{-1}(x-\beta_n) + a_n \right) P_n'.
}
$$
The matrix coefficients $\beta_n, \gamma_n, \tilde{\beta}_n, \tilde{\gamma}_n, a_n, b_n$
can be expressed in terms of $E_n$ and $\pi_n$ since
$$
\ba{l}
\beta_n = \pi_n - \pi_{n+1},
\hskip68pt
\gamma_n = E_nE_{n-1}^{-1},
\smallskip
\cr
\tilde{\beta}_{n-1} = {n-1\over n} \pi_n - {n\over n+1}\pi_{n+1},
\hskip 20 pt
\tilde{\gamma}_{n-1} = {n-1 \over n}E_n M_{n-1}M_{n-2}^{-1} E_{n-1}^{-1},
\smallskip
\cr
a_n = \beta_n - \tilde{\beta}_{n-1},
\hskip 70 pt
b_n=\gamma_n - {n \over n-1} \tilde{\gamma}_{n-1}.
\ea
$$
From here, it is just a  matter calculation to get the result, using the fact that
$M_kM_j^{-1} = \hat{M}_k\hat{M}_j^{-1} = \hat{M}_j^{-1}\hat{M}_k$, where
$\hat{M}_n=I+n\varphi_2\psi_1^{-1}.$
\epr

\subsection{Examples}

\smallskip

The purpose of the following examples is to show that non-diagonalizable matrix
$\mathcal{P}_{2,1}$-functionals do exist, even in the positive definite case, and that
the family $\mathcal{P}_{2,1}$ is strictly bigger than the zero class (excepting the
scalar case). Indeed, the presented examples are all positive definite and lie on the
class $s=1$. The matrix functionals of the examples have the structure $w(x)R(x)\,dx,$
where $w$ is a classical scalar weight and
$$
R=\pmatrix{p+qq^* & bq \cr \bar bq^* & |b|^2}, \quad p,q\in\PP,
$$
\vskip-20pt
$$
p \hbox{ with positive leading coefficient}, \quad \deg q=1, \quad b\in\C\setminus\{0\}.
$$
We will deal with a canonical form of these functionals, since any of them is congruent
to one with the form
$$
w(x)\pmatrix{\hat p(x) + |a|^2x^2 & ax \cr \bar ax & 1}dx, \quad \hat p\in\PP \hbox{
monic}, \quad a\in\C\setminus\{0\}.
$$
This kind of functionals are never diagonalizable by congruence, neither by equivalence.
This is a consequence of the fact that, as can be easily checked, any functional
$W(x)\,dx$, with
$$
W = \pmatrix{w_{11} & w_{12} \cr w_{21} & w_{22}},
$$
is non-diagonalizable by equivalence if $\{w_{11},w_{12},w_{22}\}$ is linearly
independent and $\{w_{12},w_{21}\}$ is linearly dependent.

\medskip

\noindent {\bf Example 2.} Let $u\in\PP^{{(2)}'}$ given by
$$
u = e^{-x^2} \pmatrix{1 + |a|^2x^2 & ax \cr \bar ax & 1} dx, \quad x\in\R,
\quad a\in\C\setminus\{0\}.
$$
It is not a zero class functional, but its class is $s=1$ due to the equality
$$
Du = u \pmatrix{(|a|^2-2)x & a \cr \bar a(1-|a|^2x^2) & -(|a|^2+2)x}.
$$
Besides, it is a $\mathcal{P}_{2,1}$-functional with
$\mathcal{M}_{2,1}(u) = {\rm span}_{\C^{(2,2)}}\{\Phi\}$, where
$$
\Phi(x)=\pmatrix{|a|^2+2 & 0 \cr -\bar a|a|^2x & 1}.
$$
The corresponding Pearson-type equation is $D(u\Phi)=u\Psi$, with
$$
\Psi(x)=\pmatrix{ -4x & a \cr 2\bar a & -(|a|^2+2)x}.
$$
Any right multiple of $\Phi$ by a non-singular matrix factor can be chosen as a generator
of $\mathcal{M}_{2,1}(u)$, therefore, it will play a similar role in the Pearson-type
equation for $u$. However, if we choose
$$
\Phi^{(0)}=\Phi\pmatrix{1&0 \cr 0&2},
$$
the new functional $u^{(1)}=u\Phi^{(0)}$ is again a positive definite
$\mathcal{P}_{2,1}$-functional of similar type. Indeed,
$$
u^{(1)} = e^{-x^2} \pmatrix{|a|^2 + 2 + 2|a|^2x^2 & 2ax \cr 2\bar ax & 2} dx, \quad x\in\R.
$$
This shows explicitly in the present example the general fact that any quasi-definite
$\mathcal{P}_{2,1}$-functional generates a sequence of $\mathcal{P}_{2,1}$-functionals,
according to Theorem \ref{MAIN} and Remark \ref{CHAIN}.

\medskip

\noindent {\bf Example 3.} The functional $u\in\PP^{{(2)}'}$ defined by
$$
u = x^r e^{-x} \pmatrix{x + |a|^2x^2 & ax \cr \bar ax & 1} dx, \quad x\in(0,\infty),
\quad a\in\C\setminus\{0\}, \quad r>-1,
$$
lies again on the class $s=1$ since
$$
D(uxI) = u \pmatrix{r+2+(|a|^2-1)x & a \cr -\bar a|a|^2x^2 & r+1-(|a|^2+1)x}.
$$
It is also a $\mathcal{P}_{2,1}$-functional, with
$\mathcal{M}_{2,1}(u) = {\rm span}_{\C^{(2,2)}}\{\Phi\}$ generated by
$$
\Phi(x)=\pmatrix{(|a|^2+1)x & 0 \cr -\bar a|a|^2x^2 & x}.
$$
The Pearson-type equation is $D(u\Phi)=u\Psi$, where
$$
\Psi(x)=\pmatrix{ (r+2)(|a|^2+1)-x & a \cr -(r+2)\bar a|a|^2x & r+1-(|a|^2+1)x}.
$$
Notice that $u^{(1)}=u\Phi$ is given by
$$
u^{(1)} = x^{r+1} e^{-x} \pmatrix{(|a|^2+1)x + |a|^2x^2 & ax \cr \bar ax & 1} dx,
\quad x\in(0,\infty),
$$
so, it is a positive definite $\mathcal{P}_{2,1}$-functional of similar type.

\medskip

\noindent {\bf Example 4.} The functional $u\in\PP^{{(2)}'}$ given by
$$
u = x^r e^{-x} \pmatrix{x^2 + |a|^2x^2 & ax \cr \bar ax & 1} dx, \quad x\in(0,\infty),
\quad a\in\C\setminus\{0\}, \quad r>-1,
$$
is also in the class $s=1$ since
$$
D(ux^2I) = u \pmatrix{(r+|a|^2+4)x-x^2 & a \cr -\bar a(|a|^2+1)x^2 & (r-|a|^2+2)x-x^2},
$$
and belongs to the family $\mathcal{P}_{2,1}$, with
$\mathcal{M}_{2,1}(u) = {\rm span}_{\C^{(2,2)}}\{\Phi\}$ generated by
$$
\Phi(x)=\pmatrix{x & -a \cr 0 & (r+|a|^2+2)x}.
$$
The Pearson-type equation is $D(u\Phi)=u\Psi$, with
$$
\Psi(x)=\pmatrix{ (r+|a|^2+3)-x & a \cr -\bar a(|a|^2+1)x & (r+1)(r+2)-(r+|a|^2+2)x}.
$$
As in the previous cases, there is a choice of $\Phi^{(0)}\in\mathcal{M}_{2,1}(u)$ that
makes $u^{(1)}=u\Phi^{(0)}$ a positive definite $\mathcal{P}_{2,1}$-functional of
similar type. The choice is
$$
\Phi^{(0)} = \Phi \pmatrix{r+1&0 \cr 0&1},
$$
and the new functional is then
$$
u^{(1)} = x^{r+1} e^{-x} \pmatrix{(r+1)(|a|^2+1)x^2 & (r+1)ax \cr (r+1)\bar ax & r+2} dx,
\quad x\in(0,\infty).
$$

\smallskip

\section{The zero class}

\medskip

The zero class is a specially simple subset of the family $\mathcal{P}_{2,1}.$ This
simplicity allows for zero class functionals a deeper analysis than for general
$\mathcal{P}_{2,1}$-functionals. According to the definition of the zero class we suppose
in this section that $u\in \mathbb{P}^{(m)'}$ is a quasi-definite functional that
satisfies a Pearson-type equation
\beq\label{PTEZC}
D \left(u \alpha I\right) = u \Psi,
\quad \alpha\in\mathbb{P}_2\setminus\{0\},\quad \Psi\in\mathbb{P}_1^{(m)}.
\eeq
We will use the notation $\alpha(x) = \alpha_0+\alpha_1 x + \alpha_2 x^2,$  $\alpha_i\in
\mathbb{C},$ and $\Psi(x) = \psi_0 + \psi_1x,$  $\psi_j\in\mathbb{C}^{(m,m)}.$

The first aim of this section is to obtain explicit expressions for the elements
associated with a zero class functional $u$ in terms of the coefficients
$\alpha_i\in\mathbb{C},$ $\psi_j\in\mathbb{C}^{(m,m)}.$ This will lead to a
characterization of the polynomials $\alpha \in \mathbb{P}_2\setminus\{0\},$
$\Psi\in\mathbb{P}_1^{(m)}$ which can appear in the Pearson-type equation of a zero class
functional. As a first restriction for $\alpha, \Psi,$ notice that Corollary \ref{SSC1}
implies that $\psi_1+n\alpha_2I$ must be non-singular for $n\geq 0.$

Remember that $\left(P_n\right)$ denotes the sequence of monic MOP related to $u,$
$P_n(x) = x^n I + \pi_n x^{n-1} + \cdots$ and $E_n = \langle x^n P_n, u \rangle$. As we
have shown in the proof of Corollary \ref{RRPPS}, the coefficients of the recurrence
$xP_n = P_{n+1} + \beta_n P_n + \gamma_n P_{n-1}$ and the coefficients of the relation
$P_n={1\over n+1}P_{n+1}' +a_n P_n' + b_n P_{n-1}'$ can be obtained from $\pi_n$ and
$E_n.$ So, we will just calculate $\pi_n$ and $E_n$ in terms of $\alpha$ and $\Psi.$

\medskip

From the Pearson-type equation for the functional $u$ we obtain the relation
$(\ref{RMM})$ among the moments, that can be written in the following way
\beq\label{RMO}
n\mu_{n-1} \alpha_0 + \mu_n N_n + \mu_{n+1} M_n = 0,\quad n\geq 0,
\eeq
where $N_n = \psi_0 + n \alpha_1 I,$ $M_n = \psi_1 + n \alpha_2 I.$ Taking $n=0$ and $n=1$
in $(\ref{RMO})$ we obtain
\beq\label{RMZC}
\mu_1 = -\mu_0 \psi_0 \psi_1^{-1},
\quad
\mu_2 = \mu_0\bigl(\psi_0\psi_1^{-1}\psi_0 + \alpha_1\psi_0\psi_1^{-1} - \alpha_0\bigr) M_1^{-1}.
\eeq
Let us denote $\tilde{u} = u \alpha I$ and $(\tilde{\mu}_n)_{n\geq 0}$ its corresponding
moment sequence. We know that
$$
\tilde{\mu}_n = \alpha_0\mu_n + \alpha_1\mu_{n+1} + \alpha_2\mu_{n+2}, \quad n\geq0.
$$
This equality for $n=0,$ together with $(\ref{RMZC})$, gives
$\tilde{\mu}_0 =\mu_0\alpha(-\psi_0\psi_1^{-1})\psi_1 M_1^{-1}.$
Besides, a direct calculation shows that $\pi_1=-\mu_1\mu_0^{-1}.$
So,
      $$
      \ba{l}
      \pi_1 =  E_0\psi_0\psi_1^{-1} E_0^{-1},
      \hskip 70pt
      \tilde{\pi}_n = {n\over n+1} \pi_{n+1},
      \medskip
      \cr
      \tilde{E}_0 = E_0\alpha(-\psi_0\psi_1^{-1})\psi_1 M_1^{-1},
      \hskip 30pt
      \tilde{E}_n = -{1\over n+1} E_{n+1} M_n,
      \ea
      $$
where ${1 \over n+1} P_{n+1}'(x) = x^n + \tilde{\pi}_n x^{n+1} + \cdots$ and
$\tilde{E}_n = {1\over n+1} \langle x^nP_{n+1}', \tilde{u} \rangle.$

Since $u$ is a quasi-definite  $\mathcal{P}_{2,1}$-functional, the same thing happens to
$\tilde{u}.$ Indeed, $\tilde{u}$ is also zero class because $D\left(\tilde{u}\alpha
I\right) = u\tilde{\Psi},$ $\tilde{\Psi}=\Psi + \alpha'I.$ Notice that $\tilde{\Psi}(x) =
\tilde{\psi}_0 + \tilde{\psi}_1 x$, where $\tilde{\psi}_1=M_2$ and $\tilde{\psi}_0 =N_1.$

The above results show that we can define a sequence $(u^{(j)})_{j\geq0}$ of zero class
functionals by $u^{(j)}= u \alpha^j,$ and these functionals satisfy the Pearson-type
equation
$$
D(u^{(j)} \alpha) = u^{(j)} \Psi^{(j)}, \qquad \Psi^{(j)}=\Psi + j \alpha'.
$$
Notice that $\psi_0^{(j)} = N_j,$ $ N_k^{(j)} = N_{k+j},$ $\psi_1^{(j)} = M_{2j},$
$M_k^{(j)} = M_{k+2j},$ where we denote with the superscript $(j)$ the elements
associated with the functional $u^{(j)}.$ Therefore,
        $$
        \ba{l}
        \pi_1^{(j)} = E_0^{(j)} N_j M_{2j}^{-1} (E_0^{(j)})^{-1},
        \hskip73pt
        \pi_k^{(j+1)} = {k\over k+1} \pi_{k+1}^{(j)},
        \medskip
        \cr
        E_0^{(j+1)} = E_0^{(j)} \alpha(-N_j M_{2j}^{(-1)}) M_{2j} M_{2j+1}^{-1},
        \hskip20pt
        E_k^{(j+1)} = -{1\over k+1} E_{k+1}^{(j)} M_{k+2j}.
        \ea
        $$
After an inductive process,
        $$
        \ba{l}
        \pi_n = \pi_n^{(0)} = n \pi_1^{(n-1)} = n E_0^{(n-1)} N_{n-1} M_{2n-2}^{-1} (E_0^{(n-1)})^{-1},
        \medskip
        \cr
        E_n = E_n^{(0)} = (-1)^n n! E_0^{(n)} M_{2n-2}^{-1} \cdots M_{n-1}^{-1}
            = (-1)^n n! E_0^{(n)} M_{2n-1} V_{n-1}^{-1},
        \ea
        $$
where $V_n= M_n \cdots M_{2n+1}.$ Also,
         $$
         E_0^{(n)} = E_0 \alpha(-N_0M_0^{-1})M_0 M_1^{-1} \cdots
                         \alpha(-N_{n-1}M_{2n-2}^{-1})M_{2n-2}M_{2n-1}^{-1},
         $$
         and, so,
         \beq\label{REZC}
         \kern-5pt
         E_n =(-1)^n n! E_0 \alpha(-N_0M_0^{-1})M_0 M_1^{-1} \cdots
                            \alpha(-N_{n-1}M_{2n-2}^{-1})M_{2n-2}V_{n-1}^{-1}.
         \eeq
If we define $\Pi_n = E_n^{-1}\pi_n E_n,$ then
         \beq\label{RPIZC}
         \cases{
         \Pi_n = n V_{n-1} M_{2n-2}^{-1} N_{n-1} V_{n-1}^{-1}, \medskip
         \cr
         E_n^{-1} E_{n+1} = -(n+1) V_{n-1} M_{2n-1}^{-1} \alpha(-N_nM_{2n}^{-1}) M_{2n}V_n^{-1}.}
      \eeq

The above expressions give $\pi_n$ and $E_n$ in terms of $\alpha$ and $\Psi$ for a zero
class functional $u.$ When $u$ satisfies the Pearson-type equation but it is not
quasi-definite, the expressions for $\pi_k$ and $E_k$ are valid  for the finite segment
$(P_k)_{k=0}^n$ of MOP with respect to $u,$ whenever $\Delta_0,\dots,\Delta_n$ and
$M_0,\dots,M_{2n-1}$ are non-singular. This is because, then, the previous arguments
remain valid for $(u^{(j)})_{j=0}^n$ and $(P_k^{(j)})_{k=0}^{n-j},$ as follows from
Corollary \ref{FSMOPC} and Theorem \ref{SFS}. Moreover, if $M_{2n},$ $M_{2n+1}$ are
non-singular too, the formulas are also valid for the coefficients $\pi_{n+1},$ $E_{n+1}$
of the extra polynomial $P_{n+1}$ orthogonal to $\mathbb{P}_n^{(m)},$ given by
Proposition \ref{S1P1}.

\medskip

With the achieved results we can get a characterization of the polynomials $\alpha,$
$\Psi$ related to the zero class.

\bt\label{PTEZCL}
The Pearson-type equation $D(u \alpha I) = u \Psi$, $\alpha\in\mathbb{P}_2\setminus\{0\},$
$\Psi\in\mathbb{P}_1^{(m)},$ has a quasi-definite solution $u$ if and only if $M_n$ and
$\alpha(-N_n M_{2n}^{-1})$ are non-singular for $n\geq 0,$ where $N_n=\psi_0+n\alpha_1
I,$ $M_n=\psi_1+n \alpha_2 I. $ Under these conditions, the solution of the Pearson-type
equation is unique up to left matrix factors, and the quasi-definite solutions correspond
to the non-singular choices of $\mu_0$.
\et

\bpr
If $D\left(u \alpha I\right) = u \Psi$ has a quasi-definite solution, the corresponding
matrices $E_n$ are non-singular for  $n\geq 0.$ Then, $M_n$ and $\alpha(-N_nM_{2n}^{-1})$
are non-singular for $n\geq 0,$ as can be seen from (\ref{REZC}).

For the converse, from Remark \ref{SSR}, if $M_n$ is non-singular for $n\geq 0,$ the
solutions of the Pearson-type equation are determined by the choice of $\mu_0.$ Moreover,
if, besides, $\alpha\left(-N_nM_{2n}^{-1}\right)$ is non-singular for $n \geq 0,$ the
solution $u$ is quasi-definite when $\mu_0$ is non-singular. Indeed, proceeding by
induction we can prove that there exist MOP with respect to $u$ of any degree:

\noindent $\cdot$  There exists $P_0=I,$ with $E_0=\mu_0$ non-singular.

\noindent $\cdot$  Suppose that there exists a finite segment $\left(P_k\right)_{k=0}^n$
of monic MOP with respect to $u.$ By Proposition \ref{S1P1}, there is a monic matrix
polynomial $P_{n+1}$ with $\deg P_{n+1}= n+1,$ which is orthogonal to
$\mathbb{P}_n^{(m)}.$ Since $M_k$ is non-singular for $k\geq 0$, the expression of
$E_{n+1} = \langle x^{n+1}P_{n+1}, u \rangle$ is given by (\ref{REZC}). Then, the
non-singularity of $\alpha\left(-N_k M_{2k}^{-1}\right)$ for $k \geq 0$ implies that
$E_{n+1}$ is non-singular and, hence, $\left(P_k\right)_{k=0}^{n+1}$ is also a finite
segment of MOP with respect to $u.$
\epr

\br\label{RFF}
From (\ref{REZC}), we see that the non-singularity of $M_k$ for $k \leq 2n-1$ and
$\alpha(-N_jM_{2j}^{-1})$ for $j \leq n-1,$ is equivalent to the existence of a finite
segment $\left(P_k\right)_{k=0}^n$ of MOP with respect to any solution $u$ of $D\left(u
\alpha I\right) = u \Psi$ with $\mu_0$ non-singular.
\er

\medskip

As in the classical scalar case, every matrix functional in the zero class belongs, up to
a change of variable, to one of the following types:
\vskip 0.25cm
$\bullet$  $\alpha(x) = 1 \hskip 2cm$          Hermite-type polynomials.
\vskip 0.25cm
$\bullet$  $\alpha(x) = x \hskip 2cm$          Laguerre-type polynomials.
\vskip 0.25cm
$\bullet$  $\alpha(x) = 1-x^2\hskip 1.15cm$    Jacobi-type polynomials.
\vskip 0.25cm
$\bullet$  $\alpha(x) = x^2 \hskip 1.8cm$      Bessel-type polynomials.
\medskip

\noindent The characterization given by Theorem \ref{PTEZCL} can be particularized for
any of the above canonical types. For the Hermite-type polynomials, the existence of a
sequence of MOP is equivalent to the non-singularity of $\psi_1.$ In the Laguerre case,
$\psi_1$ and $\psi_0+nI$ must be non-singular for $n\geq0.$ Jacobi-type polynomials exist
if and only if $\psi_1- nI$ and $\psi_1 \pm\psi_0 - 2nI$ are non-singular for $n\geq0,$
and, finally, the non-singularity of $\psi_0$ and $\psi_1 + nI$ for $n\geq0$
characterizes the existence of the corresponding Bessel-type polynomials. Notice that the
conditions for the existence of Hermite, Laguerre, Jacobi and Bessel-type MOP are a
natural generalization of the conditions in the scalar case.

\medskip

The non-singularity of the matrices $M_n$ appeared previously in  \cite{DG05}, as a
condition for the Hermite, Laguerre and Jacobi-type polynomials to ensure that they are
given by a Rodrigues formula. Our analysis proves that it is not necessary to impose this
condition since it is automatically satisfied by any zero class functional.

\medskip

Theorem \ref{PTEZCL} has also important practical consequences for the study of MOP. When
a matrix functional is given by a positive definite weight matrix on $\R$, the
corresponding MOP always exist. However, to decide if an arbitrary matrix of measures on
$\R$ defines a quasi-definite functional can be a hard problem, even in the hermitian
case. Theorem \ref{PTEZCL} solves this problem for any matrix functional satisfying a
Pearson-type equation like (\ref{PTEZC}). Moreover, Remark \ref{RFF} gives a
generalization that measures the length of the maximal finite segments of MOP associated
with the functional when it is not quasi-definite. Some applications of this rule can be
seen in Example 5. The importance of the above result for the zero class will be clear
later, since we will see that the only non-trivial matrix functionals in this class are
not positive definite.

\subsection{Differential equation}

\smallskip

In this section we will prove that the MOP of the zero class   satisfy a second order
differential equation that generalize the  known one in the scalar case. Notice that this
is not ensured by Theorem \ref{CHSCL} {\it (iii)}, since the right-hand side of the
differo-differential equation given by this theorem could have more than one term, as
follows from the comments in Remark \ref{S1R1}. We will also obtain the structure
relation of Theorem \ref{CHSCL} {\it (ii)}.

\medskip

In order to obtain the differential equation, starting from the study of the family
$\mathcal{P}_{2,1}$, and keeping in mind  Corollary \ref{RRPPS}, we can write for any
sequence $\left(P_n\right)$ of MOP in the zero class,
\beq\label{DE1}
P_{n\pm 1}' =  \Sigma_n^{(\pm)}P_n + \Gamma_n^{(\pm)} P_n',
\eeq
$$
\cases{
\kern-3pt \Sigma_n^{(+)} \kern-2pt = \kern-1pt (n+1) E_n M_{2n-1}^{-1} M_{n-1} E_n^{-1},
\smallskip\cr
\kern-3pt \Sigma_n^{(-)} \kern-2pt = \kern-1pt - n E_{n-1} M_{2n-1}^{-1} M_{n-2} E_n^{-1},
\smallskip\cr
\kern-3pt \Gamma_n^{(+)} \kern-2pt = \kern-1pt
(n+1) E_n M_{2n-1}^{-1} (\alpha_2 E_n^{-1} x - {1\over n} M_{2n-2} E_n^{-1} \pi_n
+ {1\over n+1} M_{2n-1} E_n^{-1} \pi_{n+1}),
\smallskip\cr
\kern-3pt \Gamma_n^{(-)} \kern-2pt = \kern-1pt
E_{n-1} M_{2n-1}^{-1} M_{n-2} E_n^{-1} (x + {1\over n} \pi_n).
\smallskip}
$$
On the other hand, Theorem \ref{CHSCL} {\it (ii)} and Remark \ref{S1R1} provide the
structure relation
\beq\label{DE2}
\alpha P_n'= n \alpha_2 P_{n+1} + \eta_n P_n + \theta_n P_{n-1},
\quad \eta_n, \theta_n\in\mathbb{C}^{(m,m)}.
\eeq
Taking derivatives in the structure relation we obtain
$$
\alpha P_n'' + \alpha' P_n' = n \alpha_2 P_{n+1}' + \eta_n P_n' + \theta_n P_{n-1}'
$$
and, using (\ref{DE1}), we get
\beq\label{DE3}
\alpha P_n'' + (\alpha'I - \Gamma_n) P_n' - \Sigma_n P_n = 0,
\eeq
$$
\cases{
\Gamma_n =  n \alpha_2 \Gamma_n^{(+)} + \theta_n \Gamma_n^{(-)} + \eta_n ,
\smallskip\cr
\Sigma_n = n\alpha_2 \Sigma_n^{(+)} + \theta_n \Sigma_n^{(-)},
\smallskip}
$$
which is the differential equation for $P_n.$

We can calculate the coefficients of the above differential equation. First of all,
notice that the coefficients $\eta_n, \theta_n$ of the structure relation can be
expressed in terms of $\pi_n$  and $E_n.$  A direct computation from the structure
relation $(\ref{DE2})$ gives
$$
\eta_n = n \alpha_1 + \left[\left(n-1\right)\pi_n - n \pi_{n+1}\right]\alpha_2,
\qquad
\theta_n=-E_nM_{n-1}E_{n-1}^{-1}.
$$
Therefore, using (\ref{DE1}), (\ref{DE3}) and the above expressions, we find
$$
\Sigma_n = n E_n M_{2n-1}^{-1} M_{n-1} \bigl[(n+1)\alpha_2 + M_{n-2}\bigr] E_n^{-1}
= n E_n M_{n-1}E_n^{-1}.
$$
In the same way, writing $\Gamma_n(x) = \Gamma_n^{(1)} x + \Gamma_n^{(0)},$
$\Gamma_n^{(i)}\in\C^{(m,m)}$, we get
$$
\Gamma_n^{(1)} = E_n M_{2n-1}^{-1} \left[n(n+1)\alpha_2^2-M_{n-1}M_{n-2}\right] E_n^{-1}
= -E_n M_{-2} E_n^{-1},
$$
$$
\kern-95pt \Gamma_n^{(0)} =
n\alpha_1 - {1\over n}E_n M_{2n-1}^{-1} \left[n(n+1)\alpha_2M_{2n-2} + \right.
$$
$$
\kern100pt + \left. M_{n-1}M_{n-2}-n(n-1)\alpha_2M_{2n-1}\right]E_n^{-1}\pi_n =
$$
$$
= n\alpha_1 - {1\over n} E_n M_{2n-2} E_n^{-1} \pi_n
= n\alpha_1 - {1\over n} E_n M_{2n-2}\Pi_n E_n^{-1},
$$
where $\Pi_n$ is given in (\ref{RPIZC}). From (\ref{RPIZC}) and the above result we
finally obtain
$$
\alpha'(x) I - \Gamma_n(x) =
E_n \psi_1 E_n^{-1} x + E_n V_{n-1} \psi_0 V_{n-1}^{-1} E_n^{-1}.
$$

\noindent Summarizing, we can enunciate the following result.

\bt \label{DifEq}
Let $u$ be a zero class functional with Pearson-type equation
$D(u\alpha) = u\Psi,$ $\alpha\in\mathbb{P}_2\setminus\{0\}, $ $\Psi\in\mathbb{P}_1^{(m)}.$

\item{(i)} If $(P_n)$ is the unique sequence of monic MOP with respect to $u,$
$$
\kern10pt
\alpha P_n'' + E_n V_{n-1} \Psi V_{n-1}^{-1} E_n^{-1} P_n' - n E_n M_{n-1} E_n^{-1} P_n = 0,
$$
\kern17pt where $M_n = \psi_1 + n \alpha_2 I$ and $V_n = M_n M_{n+1} \cdots M_{2n+1}.$

\item{(ii)} If $(Q_n)$ is the unique sequence of MOP whit respect to $u$ such that $Q_n$
${}\kern17pt$ has a leading coefficient $\kappa_n=(E_nV_{n-1})^{-1},$
$$
\kern10pt
\alpha Q_n'' + \Psi Q_n' - n M_{n-1} Q_n = 0.
$$
\et

The differential equation satisfied by the MOP of the zero class characterizes such MOP,
as the next result shows.

\bt \label{DifEq2}
Let $u$ be a zero class functional with Pearson-type equation $D\left(u \alpha I\right) =
u \Psi,$ $\alpha\in\mathbb{P}_2\setminus\{0\},$ $\Psi\in\mathbb{P}_1^{(m)}.$ Then, the
differential equation
$$
\alpha y'' + \Psi y' - n M_{n-1} y = 0
$$
has a unique (up to right matrix factors) matrix polynomic solution
$y\in\mathbb{P}^{(m)}$. This solution is the only $n$-th MOP $Q_n$ with respect to $u$
which has a leading coefficient $\kappa_n =\left(E_n V_{n-1}\right)^{-1}.$
\et

\bpr
Trying $y=\sum_{k\geq 0} c_k x^k$ as a solution of the differential equation, we obtain
the recurrence for the coefficients
$$
(n-k) M_{k+n-1} c_k = (k+1) \left[ N_k c_{k+1} + (k+2) \alpha_0 c_{k+2} \right].
$$
Since $M_n$ is non-singular for $n \geq 0,$ for every $k\neq n,$ $c_{k+1}=c_{k+2}=0$
implies $c_k=0.$ Hence, any non-trivial polynomic solution must have degree $n$, and such
a solution is determined by $c_n$.

If $c_k=0$ for $k>n$ and $c_n=\kappa_n$, there exists a unique solution that must be
$Q_n.$ If, on the contrary, $c_k=0$ for $ k>n$ but $c_n$ is arbitrary, the solution is
$Q_n L_n$, where $L_n = \kappa_n^{-1} c_n.$
\epr

\subsection{The  hermitian case}

\smallskip

Among all the zero class functionals, the hermitian ones have remarkable features that
deserve to be emphasized. Maybe one of the most important has  to do with the
diagonalizatibility.

The main purpose of this section is to prove a conjecture of Dur\'{a}n and Gr\"{u}nbaum
(see \cite{DG05}): any positive definite zero class functional is diagonalizable by congruence.
In fact, we will prove a more general result, since we will get the unitary
diagonalizability and, at the same time, under much weaker conditions for the matrix
functional. The key result to prove the referred conjecture is the following one.

\bp\label{HCDGC}
Let $u\in \mathbb{P}^{(m)'}$ be a solution of $D\left(u \alpha I\right)
= u\Psi,$ $\alpha\in\mathbb{P}_2\setminus\{0\},$ $\Psi\in\mathbb{P}_1^{(m)}.$
If $\mu_{n-2},\dots,\mu_{n+2}$ are hermitian,
$$
\psi_0^* \mu_{n+1}\psi_1 - \psi_1^* \mu_{n+1} \psi_0
= i2n(n+1)(A_0\mu_{n-1} + A_1\mu_n + A_2 \mu_{n+1}),
$$
with
$A_0=\Im(\bar{\alpha}_0\alpha_1),$
$A_1=2\Im(\bar{\alpha}_0\alpha_2),$
$A_2=\Im(\bar{\alpha}_1\alpha_2).$
\ep

\bpr
From the hypothesis,
$$
\langle \Psi^* x^n, u\Psi \rangle =  \langle \Psi^*x^n, u\Psi \rangle^*.
$$
Let us calculate $\langle \Psi^* x^n, u \Psi\rangle.$
$$
\langle \Psi^* x^n, u\Psi \rangle = \langle \Psi^* x^n, D(u\alpha)\rangle =
- n \langle \Psi^* x^{n-1}, u\alpha \rangle - \psi_1^* \langle x^n, u\alpha \rangle =
$$
$$
\kern67pt
= - n \langle \bar{\alpha} x^{n-1}, u\Psi\rangle^*
  - \bigl(\langle\bar{\alpha}x^{n-1}, u\Psi \rangle
        - \langle \bar{\alpha} x^{n-1}, u  \rangle \psi_0 \bigr)^* =
$$
$$
\kern27pt
= - (n+1) \langle \bar{\alpha} x^{n-1}, D(u\alpha) \rangle^*
  + \psi_0^* \langle x^{n-1}, u \alpha \rangle =
$$
$$
\kern18pt
= - (n+1) \langle \bar{\alpha} x^{n-1}, D(u\alpha) \rangle^*
  - {1\over n}\psi_0^* \langle x^n, u\Psi \rangle.
$$
Using the above results we get
$$
  (n+1) \bigl( \langle \bar{\alpha} x^{n-1}, D(u\alpha) \rangle
- \langle \bar{\alpha} x^{n-1},D(u \alpha) \rangle^*\bigr) =
  {1\over n} (\psi_0^* \mu_{n+1} \psi_1 - \psi_1^* \mu_{n+1}\psi_0),
$$
which, together with the equality
$$
  \langle \bar{\alpha} x^{n-1}, D(u\alpha) \rangle =
- (n-1) \langle |\alpha|^2 x^{n-2}, u \rangle
- \langle \bar{\alpha}' \alpha x^{n-1}, u \rangle,
$$
gives
$$
\psi_0^* \mu_{n+1} \psi_1 - \psi_1^* \mu_{n+1}\psi_0 =
n(n+1) \langle (\bar{\alpha}\alpha'-\bar{\alpha}'\alpha) x^{n-1}, u\rangle =
$$
$$
= i 2 n (n+1) \left[
\Im(\bar{\alpha}_0\alpha_1)\mu_{n-1} +
2\Im(\bar{\alpha}_0\alpha_2)\mu_n +
\Im(\bar{\alpha}_1\alpha_2)\mu_{n+1}
\right].
$$
\epr

Using the standard notation $[A,B]=AB-BA$ for the commutator of two square matrices $A$,
$B$, we get the following immediate consequence of Proposition \ref{HCDGC}.

\bc\label{HCDGCC}
Under the conditions of Proposition \ref{HCDGC}, if $\mu_0=I$ and $\mu_1$ is hermitian too,
$$
\psi_1^* [\mu_{n+1}, \mu_1] \psi_1 =
i 2 n(n+1) (A_0 \mu_{n-1} + A_1 \mu_n + A_2 \mu_{n+1}),
$$
with the coefficients $A_0, A_1, A_2$ as in  Proposition \ref{HCDGC}.
\ec

\smallskip

The commutativity of a set of hermitian matrices is equivalent to state that they are
simultaneously unitarily diagonalizable. Therefore, Corollary \ref{HCDGCC} relates the
possibility of diagonalizing simultaneously $\mu_n$ and $\mu_1$, to the requirement for
$\alpha$ to have real coefficients. The next theorem gives conditions which ensure that
$\alpha$ must be a real polynomial.

Remember that, if $\mu_0>0$ for a matrix functional, we can normalize it by $\mu_0=I$
without loosing any hermiticity property of the functional. So, in what follows, we will
use freely this normalization when it is possible.

\bt\label{DGC1}
Let $u\in\mathbb{P}^{(m)'}$ be a solution of $D\left(u \alpha I\right) = u \Psi,$
$\alpha\in\mathbb{P}_2\setminus\{0\},$ $\Psi\in\mathbb{P}_1^{(m)}.$ If $\mu_n = \mu_n^*$
for $n\leq 5$, then $\alpha$ is a real polynomial (up non-trivial factors) under any of
the followings conditions:

\item{(i)} $\left[\mu_2, \mu_1\right] = 0,$ $\Delta_0 > 0$ and $\Delta_1,\dots, \Delta_5$
non-singular.

\item{(ii)} $\Delta_2 > 0.$
\et

\bpr

Without loss of generality, we can suppose $\mu_0=I.$ Let $A_0, A_1, A_2$  be the
coefficients given in Proposition \ref{HCDGC}.

\item{\it (i)} $\left[E_1, \mu_1\right] =0$ since $E_1=\mu_2-\mu_1^2.$ Then, from (\ref{REZC}) for
$n=1,$ we obtain $\left[\psi_1, \mu_1\right] = 0,$ which implies
$\left[\psi_0,\mu_1\right] = 0$ because $\psi_0 = - \mu_1 \psi_1.$ Using (\ref{RMO}) and
the fact that $M_n$ is non-singular for $n\leq 3$, due to Theorem \ref{SFS}, we get
$\left[\mu_n,\mu_1\right] = 0$ for $n\leq 4.$ Then, from Corollary \ref{HCDGCC},
$$
\Delta_2 \pmatrix{A_0\cr A_1\cr A_2\cr} = 0,
$$
which implies $A_i=0,$ $\forall i.$

\item{\it (ii)} Corollary \ref{HCDGCC} for $n=1,2,3$ gives
$$
\Delta_2
\pmatrix{A_0 \cr A_1\cr A_2} =
{1\over 24 i} \pmatrix{
6 \psi_1^*[\mu_2,\mu_1] \psi_1
\cr
2\psi_1^*[\mu_3,\mu_1] \psi_1
\cr
\psi_1^*[\mu_4,\mu_1] \psi_1}.
$$
Therefore,
$$
\pmatrix{A_0 & \kern-7pt A_1 & \kern-7pt A_2} \Delta_2 \kern-2pt \pmatrix{A_0\cr A_1\cr A_2\cr}
\kern-3pt = \kern-1pt
{1\over 24i}
\psi_1^*
\left(6A_0\left[\mu_2,\mu_1\right]
+ 2A_1\left[\mu_3, \mu_1\right]
+ A_2\left[\mu_4, \mu_1\right]\right)
\psi_1.
$$
Notice that, if $P(x) = (A_0+A_1x+A_2x^2)I,$
$$
\langle P, uP^* \rangle =
\pmatrix{A_0 & \kern-7pt A_1 & \kern-7pt A_2} \Delta_2 \kern-2pt \pmatrix{A_0 \cr A_1 \cr A_2}.
$$

Let us suppose $P\neq 0.$ Since $\Delta_2>0,$ Proposition \ref{S1P3} implies that
$\langle P, uP^*\rangle>0.$ From Lemma \ref{NNSS} we know  that $\psi_1$ is non-singular,
so, the matrix $(\psi_1^{-1})^* \langle P, uP^* \rangle \psi_1^{-1}$ must be positive
definite too. On the other hand, tr$\,\left[\mu_n, \mu_1\right] = 0$ and, thus,
tr$\,((\psi_1^{-1})^* \langle P, uP^* \rangle \psi_1^{-1})=0.$ Hence, $\langle P, u
P^*\rangle$ can not be positive definite. This means that $P=0$ and $A_i=0,$ $\forall i.$
\epr

\bc\label{REAL}
For any positive definite zero class functional, the scalar polynomial of the
Pearson-type equation is real up to non-trivial factors.
\ec

The following result says that a zero class functional with a real scalar polynomial in
the Pearson-type equation does not need to many conditions to be unitarily
diagonalizable.

\bt\label{DGC2}
Let $u$ be a zero class functional with $\mu_n = \mu_n^*$ for $n\leq3$ and $\Delta_0>0.$
Then, if the scalar polynomial $\alpha$ of the Pearson-type equation is real up to
factors, $u$ is unitarily diagonalizable.

Under the above conditions, if $\mu_4, \mu_5$ are hermitian too, then, $\alpha$ is real
up to factors if and only if $u$ is unitarily diagonalizable.
\et

\bpr

Suppose without loss of generality that $\mu_0=I$. If $A_i=0,$ $\forall i,$ Corollary
\ref{HCDGCC} for $n=1$ gives $\psi_1^* \left[\mu_2,\mu_1\right] \psi_1 = 0.$ Since
$\psi_1$ is non-singular, $\left[\mu_2, \mu_1\right] = 0,$ so, there exists
$T\in\mathbb{C}^{(m,m)}$ unitary such that $T\mu_nT^*$ is diagonal for $n=1,2.$ Then,
$TE_1T^*$ is diagonal because $E_1=\mu_2-\mu_1^2.$ From (\ref{REZC}) for $n=1$ we find
that $T\psi_1T^*$ is diagonal too. Hence, $T\psi_0T^*$ is also diagonal due to the
identity $\psi_0 = -\mu_1\psi_1.$ Using (\ref{RMO}) and the non-singularity of $M_n$ for
$n\geq0$ one finds that $T\mu_nT^*$ is diagonal for $n\geq0.$

The converse when $\mu_4, \mu_5$ are hermitian follows from Theorem \ref{DGC1} {\it (i)}.
\epr

Joining Theorem \ref{DGC1} and \ref{DGC2} we achieve the following result that goes even
further than the conjecture of Dur\'{a}n and Gr\"{u}nbaum.

\bt\label{DGC3}
Let $u$ be a zero class functional with $\mu_n=\mu_n^*$ for $n\leq 5.$
Then, $u$ is unitarily diagonalizable under any of the followings conditions:

\item{(i)}  $\Delta_0 >0$ and $\left[\mu_2, \mu_1\right] =0.$

\item{(ii)} $\Delta_2 > 0.$
\et

\smallskip

Notice that some of the conditions in Theorems \ref{DGC1}, \ref{DGC2} and \ref{DGC3} can
be weakened. For example, in Theorem \ref{DGC1} {\it (i)}, it is possible to substitute
the condition $\Delta_1,\dots,\Delta_5$ non-singular by $\Delta_2$ non-singular and
$\left[\mu_3,\mu_1\right] = \left[\mu_4,\mu_1\right]=0.$

\bc \label{DGconj} (Dur\'{a}n-Gr\"{u}nbaum conjecture)
Any positive definite zero class functional is unitarily diagonalizable.
\ec

The above result does not mean that the hermitian zero class is trivial, since there
exist non-diagonalizable zero class MOP with respect to hermitian functionals which are
not positive definite (see Example 5). What is trivial is the positive definite subclass
of the zero class (actually, a bigger subclass, according to Theorem \ref{DGC3}). Hence,
positive definite Hermite, Laguerre and Jacobi-type MOP are unitarily diagonalizable.
Concerning the Bessel case we can say even something more: similarly to the scalar
situation, positive definite Bessel-type MOP do not exist, as the following proposition
asserts.

\bp\label{BESSEL}
Any zero class functional whose Pearson-type equation has a scalar polynomial with a
double root, is not positive definite.
\ep

\bpr
Assume that $u$ is a positive definite zero class functional whose corresponding
Pearson-type equation has a scalar polynomial $\alpha(x)=(x-a)^2$, $a\in\C$. From
Corollary \ref{REAL}, $a\in\R$. Also, Corollary \ref{DGconj} implies that there exists
$T\in\C^{(m,m)}$ unitary such that $T\mu_n\T^*$ is diagonal for $n\geq0$. Therefore,
$TE_1T^*$ is also diagonal and, using (\ref{REZC}), we get that $T\psi_1T^*$ and
$T\psi_0T^*$ are diagonal too. So, if we define the change of variable $t(x)=x-a$, the
diagonal hermitian matrix functional $\hat u_t=Tu_tT^*$ satisfies the Pearson-type
equation $D(\hat u_t \, t^2 I) = \hat u_t T\Psi(t+a)T^*$. Hence, $\hat u_t = \hat
u_t^{(1)} \oplus \cdots \oplus \hat u_t^{(m)}$, where $\hat u_t^{(i)}$ are scalar Bessel
functionals. Since a scalar Bessel functional can not be positive, the functional $u$ is
not positive definite, in contradiction with the hypothesis.
\epr

\subsection{Examples}

\smallskip

An example of non-diagonalizable hermitian zero class functional was presented in
\cite{CMV05}. \cite{DG05} generalizes this example and provides several non-trivial
families of hermitian matrix functionals that satisfy a Pearson-type equation like
(\ref{PTEZC}). In this section we will use these examples, including some non-hermitian
generalizations, and we will prove that the corresponding zero class MOP do exist as an
application of Theorem \ref{PTEZCL}. Notice that \cite{DG05} does not answer this
question since the analysis of the non-positive definite weights $dM(x)$ given there was
under the assumption that $\int_\R P(x)\,dM(x) P(x)$ is non-singular for any matrix
polynomial $P$ with non-singular leading coefficient, something that was not proved in
the concrete examples.

The non-diagonalizability of the functionals given in the following examples is ensured
because they have the structure $u=W(x)\,dx,$ where
$$
W = \pmatrix{w_{11} & w_{12} \cr w_{21} & 0}
$$
with $\{w_{11},w_{12}\}$ linearly independent and $\{w_{12},w_{21}\}$ linearly dependent.
These conditions imply that the functional $u$ is not diagonalizable by congruence or,
even, by equivalence.

\medskip

\noindent {\bf Example 5.} Let us consider a functional $u\in\PP^{{(2)}'}$ given by
$u=w(x)R(x)\,dx,$ where $w$ is a positive classical scalar weight with Pearson equation
$(w\alpha)'=w\beta$ and
$$
R(x) = \pmatrix{c + \int\kern-1pt{q(x)\over\alpha(x)}\,dx & a \cr b & 0},
\quad q\in\PP_1\setminus\{0\},
\quad a,b\in\C\setminus\{0\}, \quad c\in\C.
$$
Notice that $u$ is hermitian when $b = \bar a$, $c\in\R$ and $q$ is a real polynomial.

This kind of functionals always satisfy the boundary conditions which ensure that $D(u
\alpha I) = (u \alpha I)'$ (see Remark \ref{DISTR}). In fact, writing them in the
canonical representations, they have the form
$$
\ba{l}
e^{-x^2} \pmatrix{c + c_1x + c_2x^2 & a \cr b & 0} dx,
\quad x\in\R,
\smallskip
\cr
x^re^{-x} \pmatrix{c + c_1x + c_2\log(x) & a \cr b & 0} dx,
\quad x\in(0,\infty),
\smallskip
\cr
(1+x)^r(1-x)^s \pmatrix{c + c_1\log(1+x) + c_2\log(1-x) & a \cr b & 0} dx,
\quad x\in(-1,1),
\ea
$$
in the Hermite, Laguerre and Jacobi case respectively. In the above expressions
$c_1,c_2\in\C$ do not vanish simultaneously and $r,s>-1$.

The functional $u$ satisfies the Pearson-type equation
$$
D(u \alpha I) = u\Psi, \quad \Psi = \pmatrix{\beta & 0 \cr {q \over a} & \beta}.
$$
Therefore, if $q(x)=q_0+q_1x$ and $\beta(x)=\beta_0+\beta_1x$,
$$
M_n=\pmatrix{\beta_1+n\alpha_2 & 0 \cr {q_1 \over a} & \beta_1+n\alpha_2},
$$
$$
\alpha(-N_nM_{2n}^{-1}) =
\alpha\left(-{\beta_0+n\alpha_1 \over \beta_1+2n\alpha_2}\right)
\pmatrix{1 & 0 \cr \ast & 1}.
$$
Notice that, due to Theorem \ref{PTEZCL}, $\beta_1+n\alpha_2$ and
$\alpha(-{\beta_0+n\alpha_1 \over \beta_1+2n\alpha_2})$ must be different from zero for
$n\geq0$. Hence, $M_n$ and $\alpha(-N_nM_{2n}^{-1})$ are non-singular for $n\geq0$. Also,
$\mu_0$ is non-singular since
$$
\mu_0 = \nu_0\pmatrix{\ast & a \cr b & 0}, \quad \nu_0 = \int_\R w(x)\,dx.
$$
So, according to Theorem \ref{PTEZCL}, we conclude that the functional $u$ defines a
sequence of zero class MOP.

The above two-dimensional examples are only particular cases of the $m$-dimensional zero
class functionals belonging to the equivalence classes defined by
$$
\ba{l}
e^{Ax} e^{-Bx^2} dx, \quad x\in\R, \quad
\Re(\lambda)>0 \;\;\forall\lambda\in\spec(B),
\smallskip
\cr
x^A e^{-Bx} dx, \quad x\in(0,\infty), \quad
\cases{\Re(\lambda)>-1 \;\;\forall\lambda\in\spec(A),
\cr
\Re(\lambda)>0 \;\;\forall\lambda\in\spec(B),}
\smallskip
\cr
(1+x)^A (1-x)^B dx, \quad x\in(-1,1), \quad
\Re(\lambda)>-1 \;\;\forall\lambda\in\spec(A),\spec(B),
\ea
$$
where $A,B\in\C^{(m,m)}$ commute and $\spec(A)$ means the spectrum of the matrix $A$. The
restrictions for the spectra ensure the integrability for any matrix polynomial and,
together with the commutativity of $A$ and $B$, lead to a Pearson-type equation of
Hermite, Laguerre and Jacobi-type respectively, according to Remark \ref{DISTR}. The
conditions for the spectra also ensure the existence of MOP whenever $\mu_0$ is
non-singular, as follows from Theorem \ref{PTEZCL}. For some choices of $A$ and $B$ it is
possible to get an equivalent hermitian functional. This is the case of the initial
examples, as \cite{DG05} points out.

These examples do not cover the zero class functionals of Bessel-type. Such examples can
be found starting from a scalar Bessel weight. For instance, $w(x) = x^r e^{1/x},$ with
$r=-1,0,1,2,\dots,$ is a Bessel weight on the unit circle $\T:=\{x\in\C \mid |x|=1\}$
with Pearson equation $(w\alpha)'=w\beta,$ $\alpha(x)=x^2,$ $\beta(x)=(r+2)x-1.$ The
matrix function $W=wR$ satisfies the equation $(W\alpha)' = W\Psi$, where $R$ and $\Psi$
have the same meaning as previously. However,
$$
W(x) = x^r e^{1/x} \pmatrix{c + \ds{c_1 \over x} + c_2\log(x) & a \cr b & 0}
$$
is not analytic on $\T$ if $c_2\neq0.$ If, for instance, we choose a logarithm with the
discontinuity at the non-negative real axis, the matrix functional $u=W(x)\,dx,$
$x\in\T,$ verifies (see Remark \ref{DISTR})
$$
D(u \alpha I) =
(W\alpha)'(x)\,dx - i2\pi e c_2 \pmatrix{1&0 \cr 0&0} \delta(x-1)\,dx,
$$
so, it satisfies the Pearson-type equation $D(u \alpha I) = u \Psi$ when $c_2=0.$

Similarly to the initial examples, this new one is equivalent to a particular
two-dimensional case of the general $m$-dimensional zero class functionals with the form
$x^r e^{B/x} dx,$ $x\in\T,$ where $r=-1,0,1,2,\dots$ and $B\in\C^{(m,m)}$ is
non-singular. Analogously to the scalar case, these functionals satisfy a Pearson-type
equation of Bessel-type since the restriction on $r$ gives the analyticity on $\T$ for
$x^r e^{B/x}.$  As in the previous examples, the conditions for $r$ and $B$ ensure the
existence of the corresponding MOP when $\mu_0$ is non-singular, due to Theorem
\ref{PTEZCL}.

Concerning the restriction on $r$ it is known that, for the Bessel scalar case, it can be
weakened to $r \neq -2,-3,\dots$ by introducing the alternative weight on $\T$
$$
w_0(x) = \sum_{k=0}^\infty {\Gamma(r+2) \over \Gamma(r+2+k)} {1 \over x^{k+1}}.
$$
This weight satisfies the equation $(w_0\alpha)'=w_0\beta+r+1,$ $\alpha(x)=x^2,$
$\beta(x)=(r+2)x-1.$ So, according to Remark \ref{DISTR}, the scalar functional
$u_0=w_0(x)\,dx,$ $x\in\T,$ verifies the Pearson-type equation $D(u_0\alpha)=u_0\beta.$

Notice that ${\Gamma(r+2) \over \Gamma(r+2+k)} = {1 \over (r+2)_k}$ where, in general, we
denote
$$
(A)_k=\cases{I & if $k=0,$ \cr A(A+I)\cdots(A+(k-1)I) & if $k\in\N,$}
$$
for any square matrix $A.$ If $A,B\in\C^{(m,m)}$ and
$\spec(A)\cap\{0,-1,-2,\dots\}=\emptyset$, we can consider the matrix function
$$
W(x)=\sum_{k=0}^\infty (A)_k^{-1} B^k {1 \over x^{k+1}},
$$
which is analytical on $\C\setminus\{0\}.$ If, besides, $A$ and $B$ commute, then
$(W\alpha)'=W\Psi+A-I,$ $\alpha(x)=x^2,$ $\Psi(x)=Ax-B.$ Hence, the matrix functional
$u=W(x)\,dx,$ $x\in\T,$ satisfies the Pearson-type equation $D(u\alpha I) = u\Psi$
analogously to the scalar case. Therefore, Theorem \ref{PTEZCL} states that there exist
Bessel-type MOP associated with $u$ when $B$ and $\mu_0$ are non-singular.

\smallskip

\section{Other differential equations}

\medskip

Among the results proved by Dur\'{a}n in \cite{D97}, we remark in this section one concerning
the existence of differential equations for MOP with respect to hermitian functionals
$u\in\PP^{{(m)}'}$ satisfying a Pearson-type equation
$$
D(u\Phi)=u\Psi, \quad \Phi\in\PP^{(m)}_2, \quad \Psi\in\PP^{(m)}_1.
$$
The referred result states that such a Pearson-type equation, together with the
hermiticity of $u\Phi,$ is equivalent to state that the corresponding MOP $(P_n)$ satisfy
a second order differential equation
\beq \label{DE2}
P_n''\Phi^* + P_n'\Psi^* + \Lambda_n P_n = 0,
\eeq
with $\Lambda_n\in\C^{(m,m)}$ such that $\Lambda_n \langle P_n,P_n \rangle_u$ is
hermitian (actually, the result is proved in \cite{D97} for matrix orthonormal
polynomials with respect to positive definite matrix functionals, but the generalization
to the quasi-definite hermitian case is immediate). If, as in the rest of paper, we
suppose that the MOP are monic, the condition for $\Lambda_n$ becomes $\Lambda_n E_n =
E_n \Lambda_n^*.$ Also, equaling the coefficients of the highest powers of $x$ in
(\ref{DE2}) we get $\Lambda_n=-n(n-1)\psi_1^*-n\varphi_2^*=-nM_{n-1}^*.$

All the examples of $\mathcal{P}_{2,1}$-functionals $u\in\PP^{{(2)}'}$ presented in
Section 3 were hermitian and positive definite and, for all of them, we found a matrix
polynomial $\Phi\in\mathcal{M}_{2,1}(u)$ with $\det\Phi\neq0$ such that $u\Phi$ is also
hermitian and positive definite (in Examples 2 and 4 such a matrix polynomial was denoted
$\Phi^{(0)}$, we omit now the superscript for convenience). Therefore, the corresponding
MOP $(P_n)$ must satisfy a second order differential equation like (\ref{DE2}).

For instance, in the case of the functional given in Example 2
$$
u = e^{-x^2} \pmatrix{1 + |a|^2x^2 & ax \cr \bar ax & 1} dx, \quad x\in\R,
\quad a\in\C\setminus\{0\}.
$$
we find
$$
\ba{l}
P_n''(x)\pmatrix{|a|^2+2 & -a|a|^2x \cr 0 & 2} +
P_n'(x)\pmatrix{ -4x & 2a \cr 2\bar a & -2(|a|^2+2)x} +
\medskip
\cr
\kern136pt
+ \; n\pmatrix{4 & 0 \cr 0 & 2(|a|^2+2)}P_n(x)=0.
\ea
$$
This functional was previously studied in \cite{DG05b}, where it was proved that the
corresponding MOP satisfy another second order differential equation linearly independent
with respect to this one. The fact that, contrary to the scalar case, the MOP can satisfy
linearly independent second order differential equations was recently discovered (see
\cite{CG,GI}).

As for the functional
$$
u = x^r e^{-x} \pmatrix{x + |a|^2x^2 & ax \cr \bar ax & 1} dx, \quad x\in(0,\infty),
\quad a\in\C\setminus\{0\}, \quad r>-1,
$$
given in Example 3, we get
$$
\ba{l}
P_n''(x)\pmatrix{(|a|^2+1)x & -a|a|^2x^2 \cr 0 & x} +
\medskip
\cr
\kern70pt
+ \; P_n'(x)\pmatrix{(r+2)(|a|^2+1)-x & -(r+2)a|a|^2x \cr \bar a & r+1-(|a|^2+1)x} +
\medskip
\cr
\kern140pt
+ \; n\pmatrix{1 & (r+1+n)a|a|^2 \cr 0 & |a|^2+1}P_n(x)=0.
\ea
$$

Finally, Example 4 deals with the functional
$$
u = x^r e^{-x} \pmatrix{x^2 + |a|^2x^2 & ax \cr \bar ax & 1} dx, \quad x\in(0,\infty),
\quad a\in\C\setminus\{0\}, \quad r>-1,
$$
whose MOP must satisfy the differential equation
$$
\ba{l}
P_n''(x)\pmatrix{(r+1)x & 0 \cr -\bar a & (r+|a|^2+2)x} +
\medskip
\cr
+ \; P_n'(x)\pmatrix{(r+1)[(r+|a|^2+3)-x] & -(r+1)a(|a|^2+1)x \cr \bar a & (r+1)(r+2)-(r+|a|^2+2)x} +
\medskip
\cr
+ \; n\pmatrix{r+1 & (r+1)a(|a|^2+1) \cr 0 & r+|a|^2+2}P_n(x) = 0.
\ea
$$

\medskip

Let us restrict our attention now to the zero class MOP, that is, those whose
corresponding functional $u\in\PP^{{(m)}'}$ satisfies a Pearson-type equation
$$
D(u\alpha I) = u\Psi, \quad \alpha\in\PP_2\setminus\{0\}, \quad \Psi\in\PP^{(m)}_1.
$$
If $u$ is hermitian, the hermiticity of $u\alpha I$ is equivalent to saying that $\alpha$
is a real polynomial. Hence, if $u$ is hermitian and $\alpha$ is real, the MOP $(P_n)$
with respect to $u$ satisfy the second order differential equation
$$
\alpha P_n'' + P_n' \Psi^* - nM_{n-1}^* P_n = 0.
$$
This differential equation is similar, but not equal to the one given in Theorem
\ref{DifEq}. However, when $\mu_0>0$ this difference disappears since, then, Theorem
\ref{DGC2} implies that $u$ is unitarily diagonalizable. That is, there exists
$T\in\C^{(m,m)}$ unitary such that $\hat u = T u T^*$ is diagonal hermitian, so, the
corresponding monic MOP $(\hat P_n)$ must be diagonal with real polynomials in the
diagonal. Following similar arguments to those given in the proofs of the theorems in
Section 4, we find that $\hat\Psi = T \Psi T^*$ is also diagonal. Moreover, $D(\hat u
\alpha I) = u \hat \Psi$, hence, $\hat \Psi$ is real. Therefore, both differential
equations are the same for $(\hat P_n)$ and, thus, also for $(P_n)$ since $\hat P_n = T
P_n T^*.$

\medskip

Returning to the family $\mathcal{P}_{2,1},$ the two-dimensional examples that we have
found suggest that, for a big subclass of hermitian $\mathcal{P}_{2,1}$-functionals, the
related MOP satisfy a second order differential equation like (\ref{DE2}). Equivalently,
it seems that for many hermitian $\mathcal{P}_{2,1}$-functionals $u\in\PP^{{(m)}'}$ it is
possible to find a generator $\Phi$ of the module $\mathcal{M}_{2,1}(u)$ such that
$u\Phi$ is hermitian too. In particular, the referred examples seem to indicate that if
$u$ is positive definite, then $u\Phi$ is also positive definite for some generator
$\Phi$ of $\mathcal{M}_{2,1}(u).$ The characterization of the subclasses of hermitian
$\mathcal{P}_{2,1}$-functionals which are invariant under the operation $u \to u\Phi$
(for some choice of the generator $\Phi$ of $\mathcal{M}_{2,1}(u)$) remains as an open
problem. This is an important question, not only for the study of differential equations
for MOP, but also for the development of a general and systematic method to obtain
modified Rodrigues' formulas for $\mathcal{P}_{2,1}$-functionals (see \cite{DG05b} for
some examples of this kind of Rodrigues' formulas), as it will be shown in a future
paper.

\bigskip
\bigskip
\noindent{\bf Acknowledgements}
\medskip

The work of the authors was supported, in part, by a research grant from the Ministry of
Education and Science of Spain, project code MTM2005-08648-C02-01, and by Project E-64 of
Diputaci\'on General de Arag\'on (Spain).

\end{document}